\documentclass[11pt,a4paper,reqno]{amsart}
\usepackage{amsaddr}
\usepackage{amsmath,amsfonts,amssymb,amsthm,mathrsfs}
\usepackage[utf8]{inputenc}
\usepackage[T1]{fontenc}
\usepackage{lmodern}
\usepackage{latexsym}
\usepackage{cancel}
\usepackage{microtype}
\usepackage{caption}
\usepackage{MnSymbol}
\usepackage{xcolor}
\usepackage{enumerate}
\usepackage[normalem]{ulem}
\usepackage{bbm}
\usepackage{geometry}
\usepackage{marginnote}

\usepackage[colorlinks=true,linkcolor=blue,citecolor=magenta]{hyperref}

\usepackage{graphics,tikz}

\newcommand{\BB}{{\mathcal  B}}
\newcommand{\CC}{{\mathcal  C}}

\newcommand{\EE}{{\mathcal  E}}
\newcommand{\FF}{{\mathcal  F}}

\newcommand{\TT}{{\mathcal  T}}

\newcommand{\OO}{{\mathcal  O}}
\newcommand{\PP}{{\mathcal  P}}

\newcommand{\BE}{{\mathbb E}}

\newcommand{\BM}{{\mathbb M}}

\newcommand{\BP}{{\mathbb P}}
\newcommand{\BR}{{\mathbb R}}

\newcommand{\BBr}{{\mathscr B}}
\newcommand{\MMr}{{\mathscr M}}
\newcommand{\RRr}{{\mathscr R}}
\newcommand{\SSr}{{\mathscr S}}

\newcommand{\esssup}{\mathop{\mathrm{ess\,sup}}}

\newcommand{\qesssup}{\mathop{\mathrm{q\mbox{-}ess\,sup}}}

\newcommand{\fch}{{\mathbf{1}}}
\newtheorem{theorem}{\bf Theorem}[section]
\newtheorem{proposition}[theorem]{\bf Proposition}
\newtheorem{lemma}[theorem]{\bf Lemma}
\newtheorem{corollary}[theorem]{\bf Corollary}

\theoremstyle{definition}
\newtheorem{definition}[theorem]{Definition}
\newtheorem{example}[theorem]{\bf Example}

\newtheorem{remark}[theorem]{Remark}

\numberwithin{equation}{section}

\begin{document}

\title[Dirichlet problem for integro-differential equations]{Dirichlet problem for semilinear partial  integro-differential equations: the method of orthogonal projection}

\maketitle
\begin{center}
 \normalsize
  TOMASZ KLIMSIAK\footnote{T. Klimsiak: {\tt  tomas@mat.umk.pl}, A. Rozkosz: {\tt rozkosz@mat.umk.pl} (corresponding author)}\textsuperscript{1,2} \,\,\,
  ANDRZEJ ROZKOSZ \textsuperscript{2}   \par \bigskip
  \textsuperscript{\tiny 1} {\scriptsize Institute of Mathematics, Polish Academy of Sciences,\\
 \'{S}niadeckich 8,   00-656 Warsaw, Poland} \par \medskip

  \textsuperscript{\tiny 2} {\scriptsize Faculty of
Mathematics and Computer Science, Nicolaus Copernicus University,\\
Chopina 12/18, 87-100 Toru\'n, Poland }\par
\end{center}

\begin{abstract}
We  study  the Dirichlet  problem for semilinear equations on general open sets with
measure data on the right-hand side and irregular boundary data.
For this purpose we develop
the classical method of orthogonal projection.
We treat in a unified form equations with operators belonging to the
broad class of integro-differential operators associated with symmetric regular Dirichlet forms.
\end{abstract}


\section{Introduction}
\label{sec1}

\subsection{Statement of the problem}
Let $E$ be a locally compact separable  metric space, $D$ be an open subset of $E$ and $m$
be a Radon measure on $E$ with full support.
Let  $(L,\mathfrak D(L))$ be a self-adjoint operator that generates a
Markov semigroup $(T_t)_{t>0}$ on $L^2(E;m)$
and regular Dirichlet form $(\EE,\mathfrak D(\EE))$ (i.e. $L$ is a  Dirichlet operator).
The goal of  the present paper is to study, within this general framework,
the Dirichlet problem for  semilinear equations
\begin{equation}
\label{eq1.1c}
-Lu=f(\cdot,u)+\mu\quad\text{in }D.
\end{equation}
In (\ref{eq1.1c}), $f:E\times\BR\rightarrow\BR$ is a
given function and $\mu\ll\mbox{Cap}$, where  $\mbox{Cap}:2^E\to [0,\infty]$
is a  Choquet capacity associated with $L$ (note that $m\ll\mbox{Cap}$).

Let $E=\BR^d$ and $j:(0,\infty)\to [0,\infty)$ be a
Borel function. The Dirichlet problem for a special class of nonlocal (self-adjoint)
operators $L=\mathcal I$ 
that admit the form
\begin{equation}
\label{eq1.11}
\mathcal Iu(x):= \mathrm{p.v.} \int_{\mathbb{R}^d}(u(y)-u(x)) j(|x-y|)\,dy:=\lim_{\varepsilon\to 0^+} \int_{\mathbb{R}^d\setminus B(x,\varepsilon)}(u(y)-u(x)) j(|x-y|)\,dy
\end{equation}
is an area of current intensive research
(see, e.g., \cite{A,AD,BVW,BGPR,BJK,CS2,Chen,CFQ,CLO,CV,CY,FKV,K:CVPDE,KR:JFA,MSW,RS}
and the references therein).
It is by now understood that well-posed Dirichlet problem for \eqref{eq1.11}
must consist of two conditions: an exterior condition on $D^c:=\BR^d\setminus D$
and a description of the asymptotic  behavior of a solution at the boundary $\partial D$.
The last condition, in  the most general form,  has been  formulated in \cite{BVW} for a
suitable subclass of operators of the form \eqref{eq1.11}  and is based on  the notion of the
{\em boundary  trace operator} $W_D$ introduced for the fractional Laplacian  in \cite{BJK}.
For $u: D\rightarrow\bar\BR$ and $x_0\in D$ this operator is defined by
\begin{equation}
 \label{eq1.12}
W_D[u]=\lim_{V\uparrow D}\eta_V[u],\quad \eta_V[u](A)=\int_A G_V(x_0, z)\!\! \int_{D \backslash V} j(|z-y|) u(y) dy\,dz,\,\,A\in\BBr(\BR^d),
\end{equation}
whenever the measures $\eta_V[u]$ are bounded as $V\uparrow D$ and converge weakly (in the sense of measures) as $V\uparrow D$. In (\ref{eq1.12}), $G_V$ is the Green function for the operator $\mathcal I$ restricted to $V$
(in the notation we omitted the dependence of $W_D$ on $x_0$ since under the assumptions of \cite{BVW},
$G_V(x_0, z)$ and $G_V(x'_0, z)$ are suitably  comparable near $\partial V$).
The general form of  the Dirichlet condition within  the subclass of  operators  \eqref{eq1.11}  has the  form: for given functions $g$ and $h$
(or even measures) find $u$ such that
\begin{equation}
\label{eq1.1v1}
 u=g\quad\text{on } D^c,\qquad W_D(u)=h\quad\text{on }\partial D.
\end{equation}
This clearly  contrasts with the classical formulation of the Dirichlet problem in which
we are looking for a function $u:\bar D\to\mathbb R$ satisfying the equation $\Delta u=0$ in $D$  and the Dirichlet condition
\begin{equation}
\label{eq1.lap}
u=g\quad\text{on }\partial D.
\end{equation}
Our aim is to develop a unified theory of the Dirichlet problem for the
large class of Dirichlet operators in such a way that it  embraces the seemingly different boundary conditions   \eqref{eq1.1v1} (mainly with $h=0$) and \eqref{eq1.lap}. The second goal is to study this problem for semilinear equations \eqref{eq1.1c}.

The model example of a local operator satisfying our assumptions is a divergence form operator
\begin{equation}
\label{eq1.2}
L=\sum^d_{i,j=1}\partial_{x_i}(a_{ij}(x)\partial_{x_j}),
\end{equation}
where the coefficients $a_{ij}\in\mathcal B(D)$ are locally integrable, the matrix $a:=[a_{ij}]$
is nonnegative definite a.e., and either
$a$ is a.e. invertible  with $a^{-1}\in L^1_{loc}(D)$ or $a_{ij}\in W^{1,2}_{loc}(D)$, $i,j=1,\dots, d$ (see \cite{RW,VV}).
When $D$ is smooth and  $\Sigma\subset \partial D$ is a relatively open
smooth part of $\partial D$, then  the elliptic operator \eqref{eq1.2} with Neumann  boundary condition
\begin{equation}
\label{eq2.neu}
n\cdot a\nabla u+\gamma u=0\quad\text{on}\quad\Sigma
\end{equation}
also fits our framework  (see \cite{VV}, where the Wentzell
boundary problem is studied within the framework of Dirichlet forms).
As an example of purely nonlocal operator   can serve the fractional Laplacian
\begin{equation}
\label{eq1.3}
L=-(-\Delta)^{\alpha/2}
\end{equation}
with $\alpha\in(0,2)$ (see Example \ref{ex2.8}).
Another interesting example is  the so-called  regional fractional Laplacian
\begin{equation}
\label{eq1.4}
L=-(-\Delta)^{\alpha/2}_{D}
\end{equation}
(see, e.g.,  \cite{Chen} in the context of the Dirichlet problem and Example \ref{ex3}).
Of course, the class of Dirichlet operators goes far beyond the aforementioned examples.
Another examples are found for instance in  \cite{CF,FOT}. At this point we would like to mention just one  class of operators
that is covered by our theory and has recently received a lot of interest,
namely the class of  so-called mixed local and nonlocal operators, whose model example is the operator
\begin{equation}
\label{eq1.mix}
L=\Delta+\Delta^{\alpha/2}
\end{equation}
(see, e.g., \cite{BDVV2,BDVV1} and the references therein). It is worth noting
that any positive linear combination of the operators mentioned above is covered by the class considered in the paper.

In the paper, we  assume that $f$  is a Borel measurable function  such that $f(x,\cdot)$
is continuous and nonincreasing for each $x\in E$  and $f(\cdot,y)$  is quasi integrable for fixed $y\in\BR$
(each function from $L^1(E;m)$ is quasi integrable).
We would like to stress that no conditions on the growth of $f(x,\cdot)$ are imposed.
Some model examples are
\begin{equation}
\label{eq.exgen}
f(x,y)=-b(x)y|y|^{p-1},\quad f(x,y)=b(x)(1-e^y),\quad f(x,y)=b(x)(1-e^{y^2})\fch_{[0,\infty)}(y),
\end{equation}
where $p\ge 1$ and $b$ is a positive quasi integrable function
not necessarily locally integrable (for instance,  for $d\ge2$ and  $q\ge 0$ the function
$ b(x)=|x|^{-q}$
is quasi integrable for the operator \eqref{eq1.3}).
As for $\mu$, we assume that  $\mu\ll\mbox{Cap}$ and there exists
a strictly positive $\rho\in\mathbb W(D)$ 
such that $\int_D\rho\,d|\mu|<\infty$.
We use  $\mathbb W(D)$ to denote a subset  of
nonnegative  Borel measurable functions on $D$ which will be defined in Section \ref{sec.exc}.  Here  only note that it is a cone such that  $1\in\mathbb W(D)$, $u\wedge v\in\mathbb W(D)$ whenever $u,v\in\mathbb W(D)$,
and $(\mathbb W(D)-\mathbb W(D))\cap L^2(D;m)$ is dense in $L^2(D;m)$.


In order to make the exposition  of the main  results of the paper more readable, throughout the Introduction we additionally assume  that there exists the Green function $G_D$  for  the operator $L$ and domain $D$.

\subsection{Dirichlet condition.}

Since we want to cover wide class of operators
ranging from local operators (see \eqref{eq1.2}, \eqref{eq2.neu})
through a mixture of local and purely nonlocal operators (see \eqref{eq1.mix}) to purely nonlocal
operators like \eqref{eq1.3}, \eqref{eq1.4}, we must formulate  the Dirichlet condition
in a way which will unify quite extreme cases \eqref{eq1.1v1}, \eqref{eq1.lap}.
It appears that one have to consider the condition
\begin{equation}
\label{eq.dcr1}
u=g\quad\text{on}\quad \partial_{\chi}D,\qquad \hat W_D(u)=0,
\end{equation}
where $\partial_{\chi}D$, called the  harmonic boundary,  is the carrier  of the  reference harmonic measure
$\nu_m^D$ related to $L$ and $D$ (see \eqref{eq.carr} below), that  is
$\partial_{\chi}D=\{B\in\BBr(D^c): \nu_m^D(B^c)=0\}$.
By saying  "$u=g$ on $\partial_{\chi}D$" we mean $u=g$ on some element
$B\in \partial_{\chi}D$, or, equivalently, $u=g$ $\nu_m^D$-a.e. The harmonic boundary
indicates where the boundary condition sits. In different words, it indicates at which points of
$D^c$ the values of $g$ matter for the problem. This is an additional information on the
Dirichlet problem which the reader may ignore (one can  just replace $\partial_\chi D$ with
$D^c$ in all the results formulated below; however, in several places, such a replacement results in weaker assertions).
In particular, we  may take
$\partial_{\chi} D=\partial D$, $\partial D\setminus\Sigma$, $\bar D^c$, $\partial D$,
$D^c$ for the operators \eqref{eq1.2}, \eqref{eq2.neu}, \eqref{eq1.3}, \eqref{eq1.4} and \eqref{eq1.mix}, respectively.
We see that $\partial_{\chi} D$ may be equal to $\partial D$  even if $L$ is purely nonlocal.
The  operator  $\hat W_D(u)$ is an extension  of the total variation of $W_D[u]$ defined by \eqref{eq1.12}.
One may ask why in (\ref{eq.dcr1}) we do not consider more general condition $\tilde W_D(u)=h$
with a suitable extension $\tilde W_D(u)$ of \eqref{eq1.12}? The restriction to $h=0$ is the price for the unification of the Dirichlet problem.
At this point we would just like to make a remark,
postponing a more detailed discussion until introducing  some basic concepts,
that for some of the operators considered in the present paper (not only  local ones) the condition  $u=g$ on $\partial D$
completely determines the Dirichlet problem for \eqref{eq1.1c}; in a sense, for some classes of operators,
each  solution  to the Dirichlet problem for \eqref{eq1.1c} satisfies $\hat W_D(u)=0$.
Summarizing, we are looking for solutions of the problem
\begin{equation}
\label{eq1.1}
-Lu=f(\cdot,u)+\mu\quad\text{in }D,\qquad u=g\quad\text{on }\partial_{\chi}D,\qquad \hat W_D(u)=0\quad\text{on }\partial D.
\end{equation}
Of course, this is a formal expression that  we are going  to
put  into precise mathematical terms.

\subsection{Definition of a solution}

One of the goals of the present paper is to provide an analytic  definition of a solution of
\eqref{eq1.1} that ensures uniqueness.
When approaching the problem directly, it would be necessary to define the operator $\hat W_D$ and interpret appropriately the second equality in \eqref{eq1.1}, and then give right  formulation of the first equality in \eqref{eq1.1}.
We  propose another way based on the method of orthogonal projection.
At this point, however, we would like to draw attention to two substantial difficulties hidden in the problem  \eqref{eq1.1}.
One of  them lies in  the fact that we consider measure data $\mu$ on the right-hand side of the first equation in
\eqref{eq1.1}, which makes the uniqueness question for solutions of (\ref{eq1.1}) even more subtle.
Let us recall here that J. Serrin \cite{Serrin} (see also \cite{Pri}) has shown that
there exists a nontrivial   function $u\in W^{1,p}_0(D),\, p<d/(d-1)$,
such that for $L$ given by \eqref{eq1.2}
\[
-\int_D \eta\,Lu\,dm=\int_Da\nabla u\cdot\nabla \eta\,dm=0,\quad \eta\in C_c^\infty(D).
\]
This means that too small set of test functions $\eta$
in the variational approach to the first equation in \eqref{eq1.1} with $L$ given by \eqref{eq1.2} may violate  uniqueness, because it is well known that in that case we cannot expect better regularity
than $u\in W^{1,p}_0(D)$ for $p<d/(d-1)$.

The second difficulty that we would like to stress is due to the fact that within our general framework the Dirichlet operator $L$ may include local component (or may have jumps only inside $D$ as in case of \eqref{eq1.4}).
Therefore, aiming for a right definition of the Dirichlet problem  for \eqref{eq1.1c} we also have to capture
rigorously  the phenomenon that  ``$u$ reaches the boundary data $g$ at the boundary
$\partial_\chi D\cap\partial D$''. In a proper definition there should be some connection
between the values of $u$ in $D$ and the values of $g$ on $\partial_{\chi}D\cap\partial D$;
otherwise the problem is not well-posed
(a connection between the values of $u$ in $D$ and the values of $g$ on $\partial_{\chi}D\cap\bar D^c$ is provided by the operator $L$ itself).
Even in the classical case, \eqref{eq1.lap} is a rather symbolic notation for the problem of
finding a harmonic function $u:D\to\mathbb R$ that is related somehow  to the function
$g:\partial D\to \mathbb R$. If $g$ is continuous
and $D$ is  regular, then we are looking for a function $u\in C(\bar D)$
such that $u(x)=g(x)$, $x\in \partial D$, so $u$ continuously reaches the boundary   value $g$.
In general, however, where we only can  expect continuity of $u$ in $D$,
a realization of the condition ``$u=g$ on $\partial D$'' has to be suitably adjusted.
As we shall see, in the case where $L=\Delta$,
we fit into  the  theory of the Dirichlet problem
with non-regular boundary data  $g$ that is described in the monograph \cite{MV}.

Our main idea is to base the study  of \eqref{eq1.1} upon spectral synthesis.
The origins of such approach, for the Laplace operator, go back to the works of
Zaremba \cite{Zar} and Weyl \cite{We}. We show that it can be successfully
applied beyond the framework of Hilbert spaces.

Let $(\EE,\mathfrak D(\EE))$ be a regular transient Dirichlet form associated with $L$ and let $F$ denote its extended domain (then $(\EE,F)$ forms a Hilbert space). For a quasi open set $V \subset E$ (see Section \ref{sub2.1}) let $F(V)$ consists of $u\in F$ such that $u=0$ quasi everywhere (q.e.) in $V^c$, i.e. $\mbox{Cap}(V^c\cap\{u\neq 0\})=0$
(depending on the regularity of $D$, $\mbox{Cap}$ may be equivalently replaced by $m$, but not always!). $F(V)$  is a closed subspace of $F$, so the orthogonal projection operator
\[
\pi_V: F\to F(V)
\]
is well defined. At the heart of our  approach lies the fact that for
any quasi open $V\subset E$ and any  $x\in V$
there exists a bounded positive Borel measure $P_V(x,dy)$ (in fact $P_V(x,dy)\ll\mbox{Cap}$,  $x\in V$) such that
\begin{equation}
\label{eq1.beh}
P_V(u)(x):=\int_Eu(y)\,P_V(x,dy)=u(x)-\pi_V(u)(x),\quad u\in F\cap \mathcal B_b(E),\,\, m\text{-a.e. }x\in E
\end{equation}
The family $(P_V(x,dy))$ forms the so-called Poisson kernel.
This crucial result allows one to extend the operator $\pi_V$ to an operator $\Pi_V$ defined at least on $\mathcal B_b(E)\cup\mathcal B^+(E)$.
The idea of an analytic definition of a solution of  \eqref{eq1.1} is to find $u\in\mathcal B(E)$
for which there exists an increasing sequence   $(V_n)$
of  quasi open subsets of $D$ such that $\bigcup_{n\ge1} V_n=D$ q.e.   (we call it a $D$-{\em total family}) and the following conditions are satisfied:
\begin{enumerate}[(a)]
\item $P_{V_n}(|u|)<\infty$ a.e., $\Pi_{V_n}(u)\in F$, $\mathbf1_{V_n}\cdot|\mu|\in F^*$,
$\mathbf1_{V_n}\cdot|f(\cdot,u)|\in F^*$, $n\ge 1$ ($F^*$ is the dual space of $F$),
and for each $n\ge 1$,
\begin{equation}
\label{eq1.6}
-L[\Pi_{V_n}(u)]=\mathbf1_{V_n}\cdot f(\cdot,u)+\mathbf1_{V_n}\cdot\mu,\quad n\ge1,
\end{equation}
in the variational sense, i.e. for any $\eta\in F(V_n)$,
\begin{equation}
\label{eq6.2an}
\EE(\Pi_{V_n}(u),\eta)= \int_{V_n} f(\cdot,u)\eta\,dm+\int_{V_n}\eta\,d\mu,
\end{equation}

\item  $u=g$  on $\partial_\chi D$,

\item $P_{V_n}u\to P_Dg$ a.e. in $D$ as $n\to \infty$.
\end{enumerate}
Condition (c) describes what we need, i.e. that $u$ reaches $g$ at the boundary
$\partial_\chi D\cap\partial D$ and at the same time  that $\hat W_D(u)=0$ with a suitably defined $\hat W_D(u)$. 

Our basic analytic  definition of a solution is the following. We say that $u$ is
a {\em projective variational solution} of \eqref{eq1.1} if  (a)--(c) are satisfied and
$f(\cdot,u)\in L^1_\rho(D;m):=L^1(D;\rho\cdot m)$
for some  strictly positive  $\rho\in\mathbb W(D)$.

\subsection{Brief description of main results}

It appears that under natural mild assumptions on $g$ and $\mu$  projective
variational solutions coincide with probabilistic solutions defined via a Feynman--Kac
formula (see Theorem \ref{th6.7}). We will not present this result in the Introduction since
the definition of a probabilistic solution to \eqref{eq1.1} requires  probabilistic potential theory machinery (see  Section \ref{sec.exc}). However, let us mention an interesting analytic result that is a consequence of  Theorem \ref{th6.7}. Namely, under the assumption that there exists  the Green function for $L$ and $D$, the notion of probabilistic solutions is equivalent to the notion of integral solutions to \eqref{eq1.1}.
We say that $u:E\to \mathbb R$ is an {\em integral solution} of \eqref{eq1.1} if for $m$-a.e. $x\in E$,
\begin{equation}
\label{eq1.5green}
u(x)=\int_{D^c}g(y)\,P_D(x,dy)+\int_Df(y,u(y))G_D(x,y)\,m(dy)
+\int_DG_D(x,y)\,\mu(dy).
\end{equation}
(with the convention that $G_D(x,y)=0$ whenever $x\in D^c$ or $y\in D^c$).
For any Borel measure $\mu$, we let  $R^D\mu$
denote the most right term of the above equation whenever it is well defined  as
the Lebesgue integral for $m$-a.e. $x\in D$.
Our first very useful result is the following.
\medskip\\
\textbf{Theorem 1} (cf. Theorem \ref{th6.7}). {\em Assume that $P_D(|g|)<\infty$ $m$-a.e.,  $\mu\ll\mbox{\rm Cap}$,
and  $\int_D\rho\,d |\mu|<\infty$ for some strictly positive $\rho\in\mathbb W(D)$.
A function  $u\in\mathcal B(E)$ is an integral solution of \eqref{eq1.1} if and only if it is a projective variational solution of \eqref{eq1.1}.}

\medskip

The above results make it legitimate  to refer simply to solutions of  \eqref{eq1.1}
(we need not specify whether we mean  probabilistic, integral  or projective variational solutions).

The next result  says that  our solutions indeed satisfy the boundary condition formulated in \eqref{eq.dcr1}.
\medskip\\
\textbf{Theorem 2} (cf. Theorem \ref{th3.16}). {\em Let $u$ be a solution of \eqref{eq1.1}. Then, for $m$-a.e. $x\in D$,
\begin{equation}
\label{eq1.wd}
\hat W^x_D(u):= \lim_{V\uparrow D, V\subset\subset D} P_V(uR^D\kappa_D)(x)=0,
\end{equation}
where $\kappa_D$ is the  killing part of the Beurling--Deny decomposition of the form $\EE$ restricted to $D$.
Moreover, when $L=\mathcal I$ (see \eqref{eq1.11}) with $j$ satisfying the assumptions of \cite{BVW}, then
\begin{equation}
\label{eq.ident1}
\hat W^{x_0}_D(u)=W_D[u](\mathbb R^d),
\end{equation}
where $W_D[u]$ is defined by \eqref{eq1.12}.
Conversely, assume that  $u\in\mathcal B(E)$ is quasi continuous and bounded on $D\cup(\partial D\cap \partial_{\chi} D)$,
and satisfies all the conditions required in the definition of the projective variational solution of \eqref{eq1.1}
except for (c), and instead of this condition we have $\hat W^x_D(u)=0$ for $m$-a.e. $x\in D$.
Then $u$ is a solution of \eqref{eq1.1}.
}

\medskip

When $f(\cdot,u)$, $\mu \in F^*$ and $g\in F$, a right definition of \eqref{eq1.1} should agree with the usual
variational inequalities approach. In Proposition \ref{ex.prop.6.3} we show that this is true
(in that case it is enough to take $V_n=D$ in \eqref{eq1.6}).
Moreover, we provide a stability result (see Proposition \ref{prop7.2}) which
implies that each solution considered in the present paper is a limit of
variational/classical solutions (with suitable approximating sequences
$(f_n)$, $(\mu_n)$ and $(g_n)$). Let
\[
V^D(u,\eta)=2\int_{D}\int_{\mathbb R^d}(u(x)-u(y))(\eta(x)-\eta(y))j(|x-y|)\,dx\,dy.
\]
In Section \ref{sec7b}  we also show that in case $L=\mathcal I$
we have $V^D(u,u)<\infty$ provided that $V^D(P_Dg,P_Dg)<\infty$ and $f(\cdot,u),\mu \in F^*(D)$, and
\[
V^D(u,\eta)=\int_D f(\cdot,u)\eta\,dm+\int_D\eta\,d\mu,\quad \eta\in F(D).
\]
The above approach to (weak) solutions of \eqref{eq1.1} has been considered for instance
in \cite{DRV,DK,FKV,MSW}.

Our next  result concerns  uniqueness. It only requires
some of the aforementioned assumptions on the data.
\medskip\\
\textbf{Theorem 3}  (cf. Theorem \ref{prop6.6}, Corollary \ref{cor6.7}).
{\em Assume that
\begin{enumerate}
\item[\rm(A1)]  $f$  is a Borel measurable function such that $\BR\ni y\mapsto f(x,y)$ is continuous and  nonincreasing for each $x\in D$.
\end{enumerate}
Then there exists at most one solution of \eqref{eq1.1}. }

\medskip

The proof of the above result is based on the equivalence between
probabilistic and projective variational solutions of \eqref{eq1.1} stated in Theorem 1.
Thus, in fact,  it is  probabilistic. However, we also provide an analytic proof but under the additional
assumption that $\rho=1$. The advantage of the probabilistic proof is that it  is
obtained as a corollary to a much stronger result, i.e. a comparison theorem stated in Theorem \ref{prop6.6}.

The main result of the paper is an existence result for \eqref{eq1.1}.
\medskip\\
\textbf{Theorem 4} (cf. Theorem \ref{th6.11}).
{\em Let $\mu\ll\mbox{\rm Cap}$ and $\int_D\rho\,d |\mu|<\infty$ for some strictly positive $\rho\in\mathbb W(D)$.
Assume that $f$ satisfies \mbox{\rm(A1)} and
\begin{enumerate}
\item[\rm(A2)]  $f(\cdot,0)\in L^1_\rho(D;m)$ for some
 strictly positive  $\rho\in\mathbb W(D)$,
\item[\rm(A3)] $P_D(|g|)<\infty$ $m$-a.e.
\end{enumerate}
Then there exists a   solution of \eqref{eq1.1}.}

\medskip

Let us make some comments on this result.
To the best of our knowledge, except for \cite{A},
in all the previous papers on nonhomogeneous Dirichlet problem for semilinear
equations of type \eqref{eq1.1}  the following  assumption was (at least)  always made:
\begin{equation}
\label{eq1.intgh}
\int_Df(y,P_D|g|(y))G_D(x,y)\,m(dy)<\infty\quad m\text{-a.e. }x\in D.
\end{equation}
This condition seems to be natural when looking for solutions  satisfying
\begin{equation}
\label{eq1.intg}
\int_D|f(y,u(y))|G_D(x,y)\,m(dy)<\infty\quad m\text{-a.e. }x\in D
\end{equation}
(in most cases the above condition is a part  of the definition of a solution).
For example, if  $L=-(-\Delta)^{\alpha/2}$, $f(x,y)=-e^y$, $\mu=0$ and $D=B(0,1)$, then we have
\begin{equation}
\label{eq1.intst}
\int_D|f(y,u(y))|G_D(x,y)\,dy=\int_D e^{u(y)}G_D(x,y)\,dy=\int_D e^{P_Dg(y)}e^{R^Df(\cdot,u)(y)}G_D(x,y)\,dy.
\end{equation}
Observe that  $e^{R^Df(\cdot,u)}\le 1$ and we see now that  \eqref{eq1.intgh} easily implies \eqref{eq1.intg} in this case.
Condition \eqref{eq1.intgh} often appears in the literature not only to obtain \eqref{eq1.intg},
but  also plays important role in proofs in which the starting point is the
well-defined (on suitable spaces) integral operator of the form
\begin{equation}
\label{eq.gfm}
w\mapsto \int_Df(y,w(y))G_D(\cdot,y)\,m(dy).
\end{equation}
However, \eqref{eq1.intgh} is quite  restrictive. Continuing  the above example with
$g(x)=(|x|-1)^{-p}$, $|x|>1$, for some $p\in (0,1-\frac\alpha2)$, we have  by \cite[Theorem 4.2]{BVW} that
\[
P_Dg(x)\sim \delta^{-p}_D(x),\quad x\in D.
\]
Since \eqref{eq1.intgh} is equivalent to $P_Dg\in L^1_{\rho}(D)$ with $\rho=\delta_D^{\alpha/2}$ (see \cite[Theorem 1.5(iii)]{CKS}),  the above relation  implies that \eqref{eq1.intgh} does not hold. In fact, for any $q\ge 0$,
\[
\int_D |f(x,P_Dg(x))|\delta^q(x)\,dx=\infty,
\]
and even if we take $f(x,y)=-y|y|^{r-1}$, then \eqref{eq1.intgh} holds
if and only if $r<\frac{d}{p}.$

As we already mentioned, in \cite{A},  as in our paper, condition \eqref{eq1.intgh} is not assumed. Moreover,  $f$ need not be monotone.
However, in \cite{A}  only the fractional Laplace
operator is considered and it is assumed that $f$ is  continuous on  $E\times\mathbb R$,  bounded on sets of the form $E\times [a,b]$, $a\le b$,
and moreover $f\le 0$ and $g\ge 0$, $f(\cdot,0)= 0$, $\mu=0$ and $D$ is of class $C^{1,1}$.

Regarding the existence problem, one of the main results of the present paper
is the fact that we are  able to get rid of \eqref{eq1.intgh} in our general setting.
This means  in particular  that  merely  under assumptions (A1)--(A3) the competition  between
possible large  values of $u$
near the boundary  and the absorption  term $f(\cdot,u)$
always gives rise to a function $u$ which satisfies \eqref{eq1.intg}.

We prove our existence result by using methods different from those considered
before  (see, e.g., \cite{A, BVW, BJK}).
Namely, unlike  the previous papers on the subject,  we did not
try to  find a proper space for the operator \eqref{eq.gfm} and then apply a fixed point theorem.
Instead, we solve a Backward Stochastic Differential Equation (BSDE)
related to \eqref{eq1.1}. In fact,  to be precise, we provide a structure result for a
solution of this  BSDE (the existence  has been proved in \cite{K:SPA})
and then, as a corollary, we deduce the existence result for \eqref{eq1.1}.

Let us stress that in  \cite{A,BVW,BJK} the Dirichlet problem \eqref{eq1.1v1} with nonzero
$h$ is studied. This is possible due to the special form of  the operator $L$. Namely, in these papers,
\begin{equation}
\label{eq.bfo}
L=-\phi(-\Delta)
\end{equation}
for some  Bernstein function $\phi$ without drift.
We provide an existence result for a general class
of operators with Dirichlet condition \eqref{eq.dcr1}. However, for $L$  of the form
\eqref{eq.bfo} with $\phi$ satisfying some weak scaling  assumptions of  \cite{BVW},  the existence
of  a solution of  problem \eqref{eq1.1} with boundary condition \eqref{eq1.1v1} follows easily from Theorem 4.
What we want to indicate  here  is that the theory we develop  in the paper  provides an
 apparatus which
when   coupled with the knowledge of the structure  of harmonic functions with respect to the  operator $L$
yield, as straightforward conclusions, results  for the general Dirichlet problem \eqref{eq1.1v1}.
\medskip\\
\textbf{Theorem 5} (cf. Theorem \ref{th4.3}).  {\em
Assume that $L$ is of the form \eqref{eq.bfo} and $\phi$ satisfies the assumptions
of \cite{BVW}. Let $M_D:D\times \bar D\to [0,\infty)$  be the Martin kernel (see Section \ref{sec7ab})
and $\partial_mD\subset \partial D$ be the set of accessible points (see Section \ref{sec7ab}).
Let $\nu$ be a bounded Borel measure on $\partial_mD$ and $\gamma$ be a Borel measure on $D^c\setminus \partial_m D$
such that $P_D(|\gamma|)<\infty$ $m$-a.e. in $D$.
Assume that \mbox{\rm(A1)--(A3)} are satisfied and
\[
\int_D |f(y,M_D\nu(y))|G_D(x,y)\,dy<\infty,\quad x\in D.
\]
Then there exists a unique   solution to the Dirichlet problem \eqref{eq1.1c},\eqref{eq1.1v1},  with $L$ given by \eqref{eq.bfo}
and $g,h$ replaced by $\gamma, \nu$, respectively.
}

\medskip

Let $u$ be a solution of \eqref{eq1.1}.
At this point we know that there exists a $D$-total family $(V_n)$ such that
$\Pi_{V_n}(u)\in F$, $n\ge 1$.
The question is, however, what can be said about the regularity of the function $u$ itself?
This is the content of  the next  result of the paper.
\medskip\\
{\bf Theorem 6.} (cf. Theorem \ref{th5.3}). {\em Let $u$ be a solution of \eqref{eq1.1} and
$(V_n)$ be a $D$-total family such that $\Pi_{V_n}(u)\in F$, $n\ge 1$. Then
for each $n\ge 1$ and each $U\subset V_n$ such that $\mbox{\rm Cap}_{\EE^{V_n}}(U)<\infty$   (capacity of $U$ relative to $V_n$)
and $P_{V_n}(|u|)\le c$  $m$-a.e. in $U$ we have
\[
u=\eta_U\quad\mbox{on }U\quad\mbox{for some }\eta_U\in F.
\]
In particular, the above equation holds for any relatively compact nearly Borel quasi open set $U$ such that $\bar U\subset V_n$ and  $P_{V_n}(|u|)\le c$  $m$-a.e. in $U$.
}

\medskip

It is worth mentioning here that the requirement that   $P_{V_n}(|u|)\le c$ a.e. in $U$, and not necessarily on $V_n$,
is  very convenient  because for  many operators $L$   local behavior of harmonic functions ($P_{V_n}(|u|)$ is harmonic on $V_n$) is well studied in the literature.

Finally, in Section \ref{sec8} we focus  on an equivalent formulation of  the definition of a solution of \eqref{eq1.1}
that is often used in the literature for problems with  specific subclasses of  Dirichlet operators considered here and suitable data.
We show that if $\mu$ is bounded,  $P_D|g|\in L^1(E;m)$ and   $u\in L^1(E;m)$ is a solution of \eqref{eq1.1} with  $f(\cdot,u)\in L^1(D;m)$, then
 $u$ is a {\em very weak solution} of \eqref{eq1.1} with
$\mathcal C= \{\eta\in \mathfrak D(L)\cap \mathscr B_b(E): L\eta\in\mathscr B_b(E)\}$,
i.e.
\begin{equation}
\label{eq1.9}
-\int_E u\,L\eta\,dm =\int_D\eta f(\cdot,u)\,dm+\int_D\eta\,d\mu,\quad  \eta\in \mathcal C,\qquad u=g\quad\text{ on }\partial_{\chi}D
\end{equation}
(see, e.g., \cite{CV,CY}).

\subsection{Comments and related literature}

Our results are part of  the   intensively studied theory of semilinear
elliptic equations with nonlocal operators.
We deal exclusively with problems (\ref{eq1.1}) with  $\mu\ll\text{Cap}$ (there are only several  papers in the literature on nonlocal equations with true measure data)
and  $f$  nonincreasing  with respect to the second variable.
For  results for general bounded Borel measure $\mu$ but with  $g=0$
(and $f$ as in our paper) we refer the reader to \cite{CV} (with the fractional Laplacian)
and \cite{K:CVPDE} (with the same operator as here).
For the case of Laplace operator see also \cite{BB, BMP,MV} and the references therein.
As far as $f$ is concerned,  other interesting models  of type \eqref{eq1.1} with nonlocal operators
(mainly involving  fractional Laplacian) are studied in the literature  with increasing $f$  (see \cite{BCSS,BFV,CQ,FV,NS})
and nonmonotone $f$ (see \cite{BDGQ,BFRW,BDVV1,CS2,CV,QX}) but  with some natural growth  restriction
and no measure data. Finally, we  stress that the  assumption that $f(\cdot,y)$ is merely quasi integrable
for fixed $y$ permits  applying  the results of the present paper to Schr\"odinger equations (in this case $f(x,y)=-V(x)y$)
with  singular nonnegative potentials $V$ (e.g. $V(x)=1/|x|^q$ for any $q\in\mathbb R$ in the case of the fractional Laplacian).

\section{Dirichlet forms and Markov processes}
\label{sec2}

We denote by $\bar\BR$  the extended real  numbers
$\{-\infty\}\cup\BR\cup\{\infty\}$ with the usual topology.
We denote by $\BBr(E)$ the set of all Borel subsets of $E$ and by  $\BB(E)$
the set of all Borel measurable functions $u:E\to\bar\BR$.
$\BB_b(E)$ (resp. $\BB^+(E)$) is the  subset of $\BB(E)$ consisting of all bounded (resp. nonnegative) functions.  We let $\BBr^*(E)$
denote the $\sigma$-algebra of {\em universally measurable} subsets of $E$.
A set $B\subset E$ belongs to $\BBr^*(E)$ if for any probability measure $\mu$ on $\BBr(E)$
there exist $B_1,B_2\in \BBr(E)$ such that $B_1\subset B\subset B_2$ and
$\mu(B_2\setminus B_1)=0$. As usual, we set $x^+=\max\{x,0\}$, $x^{-}=\max\{-x,0\}$.

\subsection{Dirichlet forms and quasi notions}
\label{sub2.1}

We start with recalling some standards facts on Dirichlet forms. More information and details are found  for instance in \cite{CF,FOT,O,Silverstein}.
In the paper, $(\EE,\mathfrak D(\EE))$ is a regular symmetric and transient Dirichlet form on $L^2(E;m)$
and $\mathfrak D_e(\EE)$ is the extended Dirichlet space.
To simplify notation, we continue to write $F$ for $\mathfrak D_e(\EE)$. Note that $F$ with the inner product $\EE$ is a Hilbert space. The dual space of $F$ (the space of all  continuous linear functionals on $F$)
is denoted by $F^*$. We will identify a nonnegative Borel measure $\mu$ on $E$ with the linear functional on $F$ given by $u\mapsto\int_Eu\,d\mu$, whenever the integral converges for every $u\in F$.

We denote by $(L,\mathfrak D(L))$ the (unique)  self-adjoint operator on $L^2(E;m)$ corresponding to $(\EE,\mathfrak D(\EE))$. This correspondence can be characterized by
\begin{equation}
\label{eq2.7}
\mathfrak D(L)\subset \mathfrak D(\EE),\qquad \EE(u,v)=(-Lu,v),\quad u\in \mathfrak D(L),\,v\in \mathfrak D(\EE),
\end{equation}
where $(\cdot,\cdot)$ is the usual inner product in $L^2(E;m)$ (see \cite[Corollary 1.3.1]{FOT}).

In the whole paper, for a set $B\subset E$ the abbreviation q.e. in $B$  means quasi-everywhere in $B$ with respect to the capacity
$\mbox{Cap}$ associated with $\EE$ (see \cite[Section 2.1]{FOT}). If $B=E$, we write simply q.e.

Recall that a function $u$ defined q.e. on $E$ with values in $\bar\BR$ is  called {\em quasi continuous} if for any $\varepsilon>0$
there exists a closed set $F_\varepsilon\subset E$ such that $u_{|F_\varepsilon}:F_\varepsilon\to \mathbb R$ is continuous
and $\mbox{Cap}(E\setminus F_\varepsilon)<\varepsilon$.
Throughout the paper,  we always consider quasi continuous $m$-versions
of functions whenever they exist.  By \cite[Theorem 2.1.7]{FOT}, any function in $F$ admits a quasi continuous $m$-version. By \cite[Theorem 2.1.2]{FOT}, for any quasi continuous function $u$ on $E$ there exists
an increasing sequence $(F_k)$ of closed subsets of $E$ such that $u_{|F_k}$
is continuous for any $k\ge 1$  and $\mbox{Cap}(E\setminus F_k)\to 0$ as $k\to \infty$. Therefore $\mbox{Cap}(E\setminus B)=0$ with $B:= \bigcup_{k\ge 1} F_k$,
and $\mathbf1_{B}u$ is Borel measurable and quasi continuous. Consequently, each quasi continuous function has a  modification (q.e.)
that is  Borel measurable and quasi continuous.

A set $V\subset E$   is called {\em quasi open} if for any $\varepsilon>0$
there exists an open set $G_\varepsilon$ containing $V$ with $\mbox{Cap}(G_\varepsilon\setminus V)<\varepsilon$.
$\BBr^n(E)$ is the family of all {\em nearly Borel measurable} subsets of $E$ (see \cite[p. 392]{FOT}). Note that $\BBr^n(E)\subset \BBr^*(E)$.
We denote by $\mathcal O_q$ the family of all quasi open nearly Borel subsets of $E$,
and by $\mathcal O$ the family of all open subsets of $E$.  Clearly $\mathcal O\subset \mathcal O_q$.
Note that $u:E\to\bar\BR$ is quasi continuous if and only if $u$ is finite q.e.
and $u^{-1}(I)$ is a quasi open set for any open set $I\subset \mathbb R$
(see the comments preceding  \cite[Lemma 2.1.5]{FOT}).

For a given $V\in \mathcal O_q$ we denote by $(\EE^V,\mathfrak D(\EE^V))$
the Dirichlet form $(\EE,\mathfrak D(\EE))$ restricted to $V$ (see \cite[Theorem 4.4.2]{FOT}):
\begin{equation}
\label{eq.dfres}
\mathfrak D(\EE^V):=\{u\in\mathfrak D(\EE): u=0\,\text{ q.e. on }\, E\setminus V\},\quad \EE^V(u,v):=\EE(u,v),\, u,v\in\mathfrak D(\EE^V).
\end{equation}

For a quasi open $V\subset E$ we let
\[
F(V)=\{u\in F: u=0\mbox{ q.e. on }V^c:=E\setminus V\}.
\]
$F(V)$ is a closed linear subspace of $F$. We denote by $F(V)^{\bot}$ the orthogonal
complement of $F(V)$ in $F$ and by $\pi_V$ the orthogonal projection on the space $F(V)$:
\begin{equation}
\label{eq2.9}
F=F(V)\oplus F(V)^{\bot},\qquad \pi_V:F\rightarrow F(V).
\end{equation}
For $g \in F$ we set
\begin{equation}
\label{eq2.13}
h_V(g)=g-\pi_V(g).
\end{equation}
Then $h_V(g)\in F(V)^{\bot}$ and, since $h_V(g)-g\in F(V)$,
\begin{equation}
\label{eq2.14}
h_V(g)=g\quad\mbox{q.e. on }V^c.
\end{equation}

\subsection{Markov processes}
\label{sec.exc}
Let $\partial$ be a one-point compactification of $E$ if $E$ is noncompact, and  an isolated point if $E$ is compact.
We adopt the convention that every function $f$ on $E$ is extended to $E\cup\{\partial\}$ by setting $f(\partial)=0$.

We denote by $\BM=(\Omega,(\FF_t)_{t\ge0}, X=(X_t)_{t\ge0},(\theta_t)_{t\in[0,\infty]},(\BP_x)_{x\in E\cup\{\partial\}})$ a (unique) $m$-symme\-tric Hunt process with life time $\zeta$ and shift operators $\theta_t$ associated with $\EE$ in the resolvent sense, i.e. a Hunt process such that
for any $f\in \BB_b(E)\cap L^2(E)$ the resolvent of $\BM$ defined as
\begin{equation}
\label{eq2.16}
R_{\alpha}f(x)=\mathbb E_x\int^{\infty}_0e^{-\alpha t}f(X_t)\,dt,\quad x\in E,\quad \alpha>0,
\end{equation}
is a quasi continuous  $m$-version of the resolvent $G_{\alpha}f$ associated with $\EE$ (for the existence of $\BM$ see \cite[Theorem 7.2.1]{FOT}).
Here $\mathbb E_x$ denotes the expectation with respect to the measure $\BP_x$.
For  $f\in\mathcal B_b(E)$ we set
\[
P_tf(x):=\mathbb E_xf(X_t),\quad x\in E,\, t\ge 0.
\]
A universally measurable function $u:E\to [0,\infty]$ is
called $(P_t)$-excessive if $P_tu(x)\uparrow u(x)$ as $t\downarrow0$ for every $x\in E$.
By \cite[Theorem 4.6.1, Theorem A.2.7, Theorem A.2.5]{FOT}, any $m$-a.e.  finite $(P_t)$-excessive function is quasi continuous, so
it is equal q.e. to a Borel quasi continuous function.

For $f\in \BB^+(E)$ and $V\in\OO_q $ we set
\begin{equation}
\label{eq2.11}
Rf(x)=\mathbb E_x\int^{\infty}_0f(X_t)\,dt,\qquad
R^Vf(x)=\mathbb E_x\int^{\tau_V}_0f(X_t)\,dt,\quad x\in E,
\end{equation}
where
\begin{equation}
\label{eq2.23}
\tau_V=\inf\{t>0: X_t\notin V\}.
\end{equation}
By \cite[p. 392]{FOT},  $\tau_V$ is a stopping time. Note that from \cite[Theorem A.2.6, Theorem 4.1.3]{FOT} it follows that
\begin{equation}
\label{eq2.28}
\BP_x(\tau_V=0)=1 \quad \mbox{q.e. }x\in V^c.
\end{equation}
By \cite[Theorem 4.4.1]{FOT}, if $f\in\BB^+(E)$ satisfies the condition $\int_EfR^Vf\,dm<\infty$, then $R^Vf\in F(V)$ and
\begin{equation}
\label{eq2.10}
\EE(R^Vf,\eta)=(f,\eta),\quad \eta\in F(V).
\end{equation}

For a Borel signed measure $\mu$ on $E$ we denote by $|\mu|$ its total variation.
Let  $\SSr(E)$ denote the set of all {\em smooth measures} on $E$.
Recall that a nonnegative measure $\mu$ belongs to $\SSr(E)$ if  there exists
an increasing sequence $\{F_n\}$ of closed subsets of $E$ such that
$\mbox{Cap}(K\setminus F_n)\to 0$ as $n\rightarrow\infty$ for any compact $K\subset E$
and $\mathbf1_{F_n}\cdot\mu\in F^*$, $n\ge1$ (see \cite[Section 2.2]{FOT}).

Let $\mu\in\SSr(E)$. We denote by $A^{\mu}$ the unique  positive continuous additive functional of $\BM$ in the Revuz correspondence with $\mu$ (see \cite[Theorem 5.1.4]{FOT}).
For a measure $\mu$ on $E$ such that $|\mu|\in\SSr(E)$ we set $A^{\mu}=A^{\mu^+}-A^{\mu^-}$,
where $\mu=\mu^{+}-\mu^{-}$ is the Jordan decomposition of $\mu$. For $\mu\in\SSr(E)$ and
$V\in\OO_q$   we set
\begin{equation}
\label{eq2.12}
R\mu(x)=\mathbb E_xA^{\mu}_\infty,\qquad R^V\mu(x)=\mathbb E_xA^{\mu}_{\tau_V},
\quad x\in E\setminus N,
\end{equation}
where $N$ is an exceptional set for $A^{\mu}$. Since $A^{\mu}_t=\int_0^tf(X_s)\,ds$, $t\ge0$, when $\mu=f\cdot m$,
this notation is consistent with (\ref{eq2.11}). For $\mu$ such that $|\mu|\in\SSr(E)$ we set $R^V\mu=R^V\mu^+-R^V\mu^-$,
whenever this makes sense. For instance, $R^V\mu$ is well defined for q.e. $x\in E$ if $|\mu|\in\RRr(E)$. From (\ref{eq2.28}) it follows that
\begin{equation}
\label{eq2.26}
R^V\mu=0\quad\mbox{q.e. in }V^c.
\end{equation}
Note also that if $\mu\in F^*$, then $R^V\mu\in F(V)$ for $V\in\OO_q $ and
\begin{equation}
\label{eq2.20}
\EE(R^V\mu,\eta)=\langle\mu,\eta\rangle,\quad\eta\in F(V).
\end{equation}
This follows from \cite[Lemma 5.1.3, Theorem 2.2.5]{FOT} applied to the part $\EE^V$
of the form $\EE$ on $V$ and an approximation argument (see \cite[Lemma 2.2.11]{FOT}).

We set
\begin{equation}
	\label{eq2.17}
	\RRr(E)=\{\mu:|\mu|\in\mathscr S(E),R|\mu|<\infty \mbox{ q.e.}\},
\end{equation}
where $R|\mu|$ is defined by (\ref{eq2.12}).
By \cite[Lemma 2.3]{KR:CM}, in the above definition of $\RRr(E)$ one can replace q.e. by $m$-a.e.
For a positive $\rho\in\mathcal B(E)$ we denote by  $\MMr_\rho(E)$ the set of all Borel measures on $E$ such that
$\int_E\rho\,d|\mu|<\infty$. We also set $\MMr_{0,\rho}(E)=\MMr_\rho(E) \cap \SSr(E)$ and $\MMr_{0,b}(E)=\MMr_{0,1}(E)$.
In general, $\RRr(E)$ is strictly bigger than $\MMr_{0,b}(E)$ (see Examples \ref{ex7.6} and \ref{ex7.7}). In fact, by \cite[Lemma 4.10]{K:CVPDE},
\begin{equation}
	\label{eq2.17nt}
	\RRr(E)=\bigcup_\rho \MMr_{0,\rho}(E),
\end{equation}
where $\rho$ ranges over the set of all strictly positive bounded $(P_t)$-excessive functions.
It is worth noting here that a bounded signed measure $\mu$ belongs to the space $\MMr_{0,b}$ if and only if it admits the decomposition
\[
\mu=f\cdot m+\nu
\]
for some $f\in L^1(E;m)$ and $\nu\in F^*$ (see \cite{KR:BPAN}). Occasionally, for a measure $\mu$ on $E$ and a function $u$ on $E$ we will use the notation
\[
\langle\mu,u\rangle=\int_Eu(x)\,\mu(dx)
\]
whenever the integral makes sense. In the whole paper  we adopt the convention
that any Borel measure on $E$ is extended, in a standard way, to $\BBr^*(E)$.

For $V\in\mathcal O$, we denote by
$\BM^V=(\Omega,(\FF^V_t)_{t\ge0}, X=(X_t)_{t\ge0},(\mathbb P^V_x)_{x\in V\cup\{\partial\}})$ a (unique) $m$-symmetric Hunt process with life time $\zeta$ associated with $\EE^V$ in the resolvent sense (see \cite[Theorem 4.4.2]{FOT}).
We denote by  $\mathbb E^V_x$  the expectation with respect to the measure $\BP^V_x$.
$\mathbb W(V)$ is the set of all $(P^V_t)$-excessive functions.

Let $V\subset E$ and $V_n\subset E$, $n\ge1$. To simplify notation, we write $V_n\uparrow V$ q.e.
if $V_n\subset V_{n+1}\subset V$ for $n\ge1$ and $\bigcup_{n\ge1}V_n=V$ q.e.

\begin{lemma}
\label{lem2.1}
Suppose that $\mu\in\SSr(E)$. Then there exists $(V_n) \subset\OO_q $ such that $V_n\uparrow E$ q.e.
and $\mathbf1_{V_n}\cdot\mu\in F^*$,
$R(\mathbf1_{V_n}\cdot\mu)\in \mathcal B^n_b(E)$, $n\ge 1$.
\end{lemma}
\begin{proof}
Let $f\in\mathcal B_b^+(E)\cap L^1(E;m)$ be a strictly positive function such that
$Rf\le1$. Set
\[
\varphi(x)= \mathbb E_x\int_0^\infty e^{-A^\mu_t}f(X_t)\,dt,\quad x\in E.
\]
By \cite[Lemma 5.1.5(ii)]{FOT}, $\varphi$ is quasi continuous and
\begin{equation}
\label{eq2.24}
\mathbb E_x\int_0^\infty \varphi(X_t)\,dA^\mu_t\le Rf(x),\quad x\in E.
\end{equation}
Moreover, by \cite[Theorem 4.6.1]{FOT}, there is a set $N$ such that $\mbox{Cap}(N)=0$ and $\varphi|_{E\setminus N}$
is nearly Borel measurable and quasi continuous. Let $V_n=\{x\in E\setminus N: \varphi(x)>1/n\}$. Then $V_n\in\OO_q$.
Observe that $V_n\uparrow \{x\in E\setminus N:\varphi(x)>0\}=E\setminus N$.
Thus  $\bigcup_{n\ge 1}V_n=E$ q.e.  By \eqref{eq2.24},
\[
R(\mathbf 1_{V_n}\cdot\mu)\le nR(\varphi\cdot\mu)\le nRf.
\]
From this and the assumptions made on $f$ one easily deduces
that $\mathbf1_{V_n}\cdot\mu\in F^*$.
\end{proof}

Let $\PP$ denote the set of all probability measures on $\BBr(E)$ and let $\FF^0_{\infty}=\sigma(X_t,t\ge0)$.  For $\mu\in\PP$
we set
\[
\BP_\mu(\Lambda)= \int_E\BP_x(\Lambda)\,\mu(dx),\quad \Lambda\in\FF^0_{\infty}.
\]
The expectation with respect to $\BP_{\mu}$ will be denoted by  $\mathbb E_{\mu}$.

\begin{definition}
We say that a family $\{P(x,dy),\, x\in E\}$ is a {\em sub-stochastic kernel }  if
\begin{enumerate}[{\rm (a)}]
\item $E\ni x\mapsto P(x,B)$ is universally measurable for any $B\in\BBr(E)$,
\item for each $x\in E$, $\BBr(E)\ni B\mapsto P(x,B)$ is a  smooth measure with $P(x,E)\le 1$.
\end{enumerate}
\end{definition}

By \cite[Theorem 4.3.2]{FOT}  (see also \cite{Silverstein}),
for any $V\in\mathcal O_q$ the  family
\begin{equation}
\label{eq2.8}
P_V(x,B):=\BP_x(X_{\tau_V}\in B),\quad x\in E, \, B\in\BBr(E),
\end{equation}
defines a sub-stochastic kernel  and for any  $g\in F$,
\begin{equation}
\label{eq2.2}
h_V(g)(x)=\int_{V^c}g(y)\,P_V(x,dy)\quad \mbox{q.e. }x\in E.
\end{equation}
For each  $g\in\BB^+(E)$ (or  $g\in\BB_b(E)$) we let
\begin{equation}
\label{eq2.3}
P_V(g)(x)= \int_{V^c} g(y)\,P_V(x,dy),\quad x\in E.
\end{equation}
Observe that
\begin{equation}
\label{eq2.25}
P_V(g)(x)=\mathbb E_x g(X_{\tau_V}),\quad x\in E.
\end{equation}
For  $g\in\BB^+(E)$ we  let
\begin{equation}
\label{eq2.19}
\Pi_V(g)(x) =g(x)-P_V(g)(x),\quad x\in E.
\end{equation}
Clearly, $P_V(x,dy)$ is concentrated on $V^c$, but if $x\in V$ and $X$ has continuous sample paths,
i.e. when $\EE$ is local (see \cite[Theorem 4.5.1]{FOT}), it is concentrated
on the topological boundary $\partial V$.  Note also that by (\ref{eq2.28}), for any $g\in \BB^+(E)$ and $V\in\OO_q$\,,
\begin{equation}
\label{eq2.8reg}
P_V(g)(x)=g(x)\quad \mbox{q.e. }x\in E\setminus V.
\end{equation}

For $W\in \mathcal O_q$ we set
\begin{equation}
\label{eq2.18}
\RRr(W)=\{\mu:|\mu|\in\SSr(E), R^W|\mu|<\infty\mbox{ q.e.}\}.
\end{equation}
Since $R^E=R$, this notation is consistent with (\ref{eq2.17}). Elements of $\RRr(W)$ may be called smooth (signed) measures of finite potential on $W$.
By \cite[Proposition 3.2]{KR:CM} applied to the form $\EE^W$, $\MMr_{0,b}(W)\subset\RRr(W)$.

The following two simple lemmas will be useful.

\begin{lemma}
\label{lem2.4}
Let $V,W\in\mathcal O_q $ and $V\subset W$. If $\mu\in\RRr(W)$, then
$\Pi_V(R^W\mu)=R^V\mu$ q.e.
\end{lemma}
\begin{proof}
Without loss of generality (see Lemma \ref{lem2.1}), we may assume that $\mu\ge0$ and $\mu\in F^*$.
Let $\eta\in F(V)$. Then
\[
\EE(R^W\mu,\eta)=\int_V\eta\,d\mu=\EE(R^V\mu,\eta).
\]
Hence $\EE(R^W\mu-R^V\mu,\eta)=0$ for $\eta\in F(V)$, which implies that $\Pi_V(R^W\mu-R^V\mu)=0$ q.e. As a result,   $\Pi_V(R^W\mu)=R^V\mu$ q.e.
\end{proof}
Note that Lemma \ref{lem2.4} is a slight generalization of Dynkin's formula (see \cite[(4.4.3)]{FOT}).

\begin{lemma}
\label{lem2.5}
Let $g\in F$. If $V,W\in\mathcal O_q $ and $V\subset W$, then
$P_V(P_W(g))=P_W(g)$ q.e.
\end{lemma}
\begin{proof}
Set $w=\Pi_V(P_W(g))$. Since $\Pi_V$ is a self-adjoint (as a projection)
operator and $w\in F(V)\subset F(W)$, $P_W(g)\in F(W)^{\bot}$, we have
\[
\EE(w,w)=\EE(P_W(g),\Pi_V(P_W(g)))=0,
\]
which implies the  desired result.
\end{proof}

\begin{corollary}
\label{cor2.6}
For any  $V,W\in\mathcal O_q $ such that $V\subset W$ we have
\[
\int_EP_V(x,dz)P_{W}(z,dy)=P_W(x,dy)\quad\text{for q.e. }x\in E.
\]
\end{corollary}
\begin{proof}
Set $\mu_x(dy)=\int_EP_V(x,dz)P_{W}(z,dy)$ and $\nu_x(dy)=P_W(x,dy)$.
By Lemma \ref{lem2.5}, for any $f\in C_c(E)\cap F$, $\langle\mu_x,f\rangle=\langle\nu_x,f\rangle$ for q.e. $x\in E$
(we use  separability of $C_c(E)$).
Since $(\EE,\mathfrak D(\EE))$ is regular, using an approximation argument we get
the above equality for all $f\in C_c(E)$.  This implies the desired result.
\end{proof}

For   $V\in \mathcal O_q$ we define a Borel measure $\nu^V_m$  on $E$  by
\begin{equation}
\label{eq.carr}
\nu^V_m(A)= \int_V P_V(x,A)\,m(dx).
\end{equation}
We call $\partial_{\chi}V=\{B\in\BBr(V^c):\nu^V_m(B^c)=0\}$ the harmonic boundary of $V$.
By writing  "$u=g$ on $\partial_{\chi}V$" we mean that $u=g$ on some element   $B\in \partial_{\chi}V$.

\begin{example}
\label{ex2.7}
Let $E:=\BR^d$, $d\ge3$, and $m$ be the Lebesgue measure on $\BR^d$.\\
(i) (Divergence form operator). Let  $a_{ij}:\BR^d\rightarrow\BR$ be measurable functions such that
\[
a_{ij}=a_{ji},\qquad \Lambda^{-1}|\xi|^2\le \sum^{d}_{i,j=1}a_{ij}(x)\xi_i\xi_j\le\Lambda|\xi|^2,\quad x,\xi \in\BR^d,
\]
for some $\Lambda\ge1$. Consider the Dirichlet form on $L^2(\BR^d;m)$ defined by
\[
\EE(u,v)=\sum^d_{i,j=1}\int_Da_{ij}(x)\frac{\partial u}{\partial x_i}(x)
\frac{\partial v}{\partial x_j}(x)\,dx,\quad u,v\in \mathfrak D(\EE):=H^1(\BR^d),
\]
where $H^1(\BR^d)$ is the usual Sobolev space of order 1. It is regular and transient (see \cite[Example 1.5.2]{FOT}).
The operator associated with $(\EE,\mathfrak D(\EE))$ in the sense of (\ref{eq2.7})
is formally  given by (\ref{eq1.2}) (this is one of possible definitions of \eqref{eq1.2}).
By \cite[Example 1.5.2]{FOT}, the extended space can be characterized as follows:
\[
F=H^1_e(\BR^d):=\{u=N*f:f\in L^2(\BR^d)\},
\]
where $N(x)=c_{d}|x|^{2-d}$ is the Newtonian kernel ($c_{d}$ is a positive constant).
By \cite{CFMS},  $P_V(x,dy)$ may by singular with respect  to the surface measure on $\partial D$ (even if $D$ is smooth),
but always, if $D$ is connected, then the support of $P_V(x,dy)$ equals $\partial D$ for $x\in D$.
By \cite[Example 2.3.2]{FOT}, for every $V\in\OO_q$,
\[
F(V)=H^1_{e,V}:=\{u\in H^1_e(\BR^d):u=0\mbox{ q.e. on }\BR^d\setminus V\}=H^1_{0,e}(V),
\]
where $H^1_{0,e}(V)$ is the extended space of $H^1_0(V)$. If $V$ is bounded, then by
Poincar\'e's inequality, $H^1_{0,e}(V)=H^1_0(V)$. For a different characterization of $F$ we refer the reader to \cite[Example 1.5.3]{FOT}.

The process $\BM$ associated with $\EE$ is a conservative diffusion (see \cite[Example 4.5.2]{FOT}).
In particular, if $V$ is open and bounded, then $\BP_x(\tau_V<\infty)=1$ and $X_{\tau_V}\in\partial V$
for every $x\in V$ (we can take $N=\emptyset$ in (\ref{eq2.2})--(\ref{eq2.8})). From (\ref{eq2.8})
it follows that  $P_V(x,dy)$ is the harmonic measure on the topological boundary $\partial V$ associated with $L$.
\smallskip\\
(ii) (Laplace operator). Consider now the special case where  $a_{ij}=\delta^i_j$, i.e. $L=\Delta$.
Then  $\BM$ is a Brownian motion running with a time clock twice as fast as the standard one.
If $V$ is open and regular, say of class $C^2$, the harmonic measure  has a strictly  positive
density $P_V(x,y)$ with respect to the surface measure $\sigma$ on $\partial V$ (see, e.g., \cite[Section 1.4]{CZ}):
\[
P_V(g)=\int_{\partial V}g(y)P_V(x,y)\,\sigma(dy),\quad x\in D.
\]
\end{example}

\begin{example}[Fractional Laplacian]
\label{ex2.8}
Let $\alpha\in(0,2)$ and $m$ be the Lebesgue measure on $\BR^d$, $d>\alpha$.
Consider the form  on $L^2(\BR^d;m)$ defined by
\[
\begin{cases}
\EE(u,v)=\int_{\BR^d}\hat u(x)\overline{\hat v(x)}
|x|^{\alpha}\,dx,\quad u,v\in \mathfrak D(\EE):=H^{\alpha/2}(\BR^d),
\smallskip\\
H^{\alpha/2}(\BR^d)=\{u\in L^2(\BR^d):\int_{\BR^d}|\hat
u(x)|^2|x|^{\alpha}\,dx<\infty\}
\end{cases}
\]
($\hat u$ is the Fourier transform of $u$). Equivalently, $(\EE,\mathfrak D(\EE))$ can be
defined by (\ref{eq2.21}), (\ref{eq2.22}) below with $D:=\BR^d$. It is a transient symmetric
regular Dirichlet form (see \cite[Example 1.4.1, Example 2.3.1]{FOT}).
The operator associated with it in the sense of (\ref{eq2.7}) is the fractional Laplace operator $-(-\Delta)^{\alpha/2}$.
By \cite[Example 1.5.2]{FOT}, the extended space can be characterized as follows:
\[
F=H^{\alpha/2}_e(\BR^d):=\{u=I_{\alpha}*f:f\in L^2(\BR^d)\},
\]
where $I_{\alpha}(x)=c_{d,\alpha}|x|^{\alpha-d}$ is the Riesz convolution kernel ($c_{d,\alpha}$ is a positive constant). Hence
\[
F(V)=H^{\alpha/2}_{e,V}:=\{u\in H^{\alpha/2}_e(\BR^d):u=0\mbox{ q.e. on }\BR^d\setminus V\}.
\]
The process associated with $\EE$ is the rotation invariant  $\alpha$-stable process.
For $x\in \BR^d$ the distribution $P_V(x,dy)$  is called the $\alpha$-harmonic measure.
It is concentrated on $V^c$. If $V$ is open, then for $x\in V$ it is absolutely continuous in
the interior of $V^c$ with respect to the Lebesgue measure. Its density function $P_V(x,y)$,
called the Poisson kernel, is strictly positive (see, e.g., \cite[(4.13)]{BRSW}.
If $V$ has the outer cone property, then $P_V(x,\partial V)=0$ for $x\in V$
(see \cite[Lemma 6]{B}). In particular, $P_V(x,dy)=P_V(x,y)\,dy$ on the whole of $V^c$.
By the aforementioned properties of $P_V$ we have  $\partial_{\chi}V=\bar V^c$ $\nu^V_m$-a.e.
for $V$ having the outer cone property.
\end{example}

\begin{example}[Regional fractional Laplacian]
\label{ex3}
Let $\alpha\in(0,2)$, $m$ denote the $d$-dimensional Lebesgue measure and $D\subset \BR^d$ be a
$d$-set, i.e. there exist  constants $c_1,c_2>0$
such that for any $x\in D$ and $r\in (0,1]$ we have $c_1r^d\le m(B(x,r)\cap D)\le c_2r^d$. Consider the form
\begin{equation}
\label{eq2.21}
\mathcal E(u,v)=c\int_D\int_D\frac{(u(x)-u(y))(v(x)-v(y))}{|x-y|^{d+\alpha}}\,dx\,dy,
\quad u,v\in \mathfrak D(\mathcal E),
\end{equation}
\begin{equation}
\label{eq2.22}
\mathfrak D(\mathcal E)= \Big\{u\in L^2(D;m): \int_D\int_D\frac{(u(x)-u(y))^2}{|x-y|^{d+\alpha}}\,dx\,dy<\infty\Big\}.
\end{equation}
It is a regular Dirichlet form on $L^2(\bar D;m)$ (see \cite[Theorem 2.2, Remark 2.1]{BBC}).
The  operator $(-L,\mathfrak D(L))$ associated with $(\EE,\mathfrak D(\EE))$
is called the regional fractional Laplacian.
\end{example}

\begin{remark}
Let $\lambda>0$ and $L_{\lambda}=L-\lambda$, where $L$ is defined by (\ref{eq2.7})
for some symmetric and  regular but  not necessarily transient
Dirichlet form $\EE$ on $L^2(E;m)$.  Then $L_{\lambda}$ corresponds,
in  the sense of (\ref{eq2.7}), to the symmetric,   regular, transient Dirichlet form
$(\mathcal{E}_{\lambda},\mathfrak D(\EE))$, where
\begin{equation}
\label{eq2.15}
\mathcal{E}_{\lambda}(u,v)=\EE(u,v)+\lambda(u,v),\quad u,v\in \mathfrak D(\EE).
\end{equation}
Therefore the results of the paper apply to the operator $L_{\lambda}$
for any symmetric regular Dirichlet form $\EE$.
Note also that the extended Dirichlet space of $(\mathcal{E}_{\lambda},\mathfrak D(\EE))$
coincides with $\mathfrak D(\EE)$ (see \cite[Theorem 1.5.3]{FOT}). Consequently, in the
case where problem (\ref{eq1.1}) with $L$ replaced by $L_{\lambda}$ is considered, we have $F=\mathfrak D(\EE)$.
\end{remark}

\section{Dirichlet problem for semilinear equations}
\label{sec6}

Throughout this section,  $D$ is an open (nonempty) subset of $E$.
We assume as given Borel measurable functions $f:E\times\BR\rightarrow\BR$,
$g:E\rightarrow\BR$ and a smooth measure $\mu$ on $D$.
Our aim is to show  an existence and   uniqueness  result for  problem \eqref{eq1.1}
in the case where  the data satisfy hypotheses (H1)--(H4) given below.

A function $h\in\BB(D)$ is said to be quasi integrable ($h\in qL^1(D;m)$ in abbreviation)
if $\BP_x(\int_0^{\tau_D}|h(X_t)|\,dt<\infty)=1$ for q.e. $x\in D$. Note that
\begin{equation}
\label{eq6.6}
L^1(D;m)\subset qL^1(D;m)\subset \SSr(D)-\SSr(D).
\end{equation}
In fact, if $h\cdot m\in\RRr(D)$ ($L^1(D;m)\subset \RRr(D)$, see (\ref{eq2.18})), then $h\in qL^1(D;m)$.
Indeed, $t\mapsto\int_0^{t\wedge\tau_D}|h(X_s)|\,ds$ is a positive continuous additive functional of $\BM^D$
in the Revuz correspondence with the measure $|h|\cdot m$. Hence, if  $h\cdot m \in \RRr(D)$, then $h\in qL^1(E;m)$
by the definition of $\RRr(D)$. The second inclusion in (\ref{eq6.6}) is an immediate consequence of the Revuz
correspondence (see \cite[Theorem 5.1.4]{FOT}).

Our basic assumptions on $f,g$ are the following.
\begin{enumerate}
\item[(H1)]  $\BR\ni y\mapsto f(x,y)$ is continuous and  nonincreasing for each $x\in D$.

\item[(H2)] $f(\cdot,y)\in qL^1(D;m)$ for each $y\in \BR$ and $f(\cdot,0)\cdot m\in\RRr(D)$.

\item[(H3)]$g\in\BB(E)$ and $P_D(|g|)<\infty$ $m$-a.e. (equivalently, q.e.).
\item[(H4)] $\mu\in\RRr(D)$.
\end{enumerate}

It is worth noting that if (H2), (H4) are satisfied, then $f(\cdot,0)\cdot m+\mu\in\RRr(D)$, and if $f$ satisfies (H1), then $f(\cdot,\cdot)-f(\cdot,0)$ satisfies (H1) as well. Therefore in the study of  (\ref{eq1.1}) (under the above assumptions) one can assume without loss of generality that $f(\cdot,0)=0$.

\subsection{Projective variational and probabilistic solutions.}

\begin{definition}
\label{def5.1}
Let $W\in\OO_q $. We say that a  family $\mathcal S=(V_n)_{n\ge1}\subset \mathcal O_q $ is {\em $W$-total} if
$V_n\subset V_{n+1}$ for $n\ge1$ and
$\bigcup_{n\ge 1}V_n=W$ q.e.
\end{definition}

The underlying definition of a solution to \eqref{eq1.1} shall  be the one below
based on the spectral synthesis.

\begin{definition}[Projective variational solutions]
\label{def6.1}
We say that  $u\in\mathcal B(E)$  is a  solution of \eqref{eq1.1} if
\begin{enumerate}[(a)]
\item $f(\cdot,u)\in L^1_\rho(D;m)$ for some strictly positive $\rho\in\mathbb W(D)$ (see \eqref{eq2.17nt}) and there exists a $D$-total family  $\mathcal{S}=(V_n)$
such that for each $n\ge 1$ we have $P_{V_n}(|u|)<\infty$ $m$-a.e., $\Pi_{V_n}(u)\in F$,
$\mathbf1_{V_n}|f(\cdot,u)|\in F^*$, $\mathbf1_{V_n}\cdot |\mu|\in F^*$, and moreover,
 for every $\eta\in F(V_n)$,
\begin{equation}
\label{eq6.2}
\EE(\Pi_{V_n}(u),\eta)= \int_{V_n} f(\cdot,u)\eta\,dm
+\int_{V_n}\eta\,d\mu,
\end{equation}

\item $u=g$  on $\partial_{\chi} D$,

\item $P_{U_n}(u)\to P_{D}(g)$ q.e.  in $D$ for any $D$-total family  $(U_n)$ satisfying (a).
\end{enumerate}
\end{definition}

\begin{remark}
In general, $\mathcal S\subset\mathcal O_q$, which  of course  does not exclude
the situation where $\mathcal S\subset\mathcal O$.  In the latter case to apply the theory
presented in the present paper it is enough to know that there
exists a family of kernels \eqref{eq2.2} but only for $V\in\mathcal O$.
The construction of such family follows in an elementary way
from Riesz's theorem (see \cite[Theorem 6.19, p. 130]{Rudin}) and Proposition \ref{prop3.2}
provided that we know that bounded  harmonic functions on open sets are continuous.
Indeed, it is enough to observe that   $h_Vg$ is a harmonic function on $V\in\mathcal O$ for any $g\in F$, i.e.   $\EE(h_Vg,\eta)=0,\, \eta \in F(D)$. Then, under the aforementioned assumptions, we have
$h_V: C_0(E)\cap F\to C_0(E)$. Applying now  Proposition \ref{prop3.2} and Riesz's theorem
yields the existence of a family of kernels $\{P_V(x,dy),\, x\in E,\, V\in\mathcal O\}$ satisfying \eqref{eq2.2}
(even for every $x\in E$).
\end{remark}

\begin{remark}
\label{rem3.ext}
If $\mathcal S$ is a total family appearing in condition (a) of Definition \ref{def6.1}
and $V\in \mathcal S$, then (a) holds for any $U\in\mathcal O_q$ such that $U\subset V$.
Indeed, the fact that   $\mathbf1_U|f(\cdot,u)|$, $\mathbf1_U\cdot|\mu|\in F^*$  is trivial.
Next, we have $\Pi_V(u)\in F$. Hence  $|\Pi_V(u)|\in F$ and
\[
|\Pi_V(u)|=|u-P_V(u)|\ge |u|-P_V(|u|).
\]
Therefore,  by   Lemma \ref{lem2.5},
\[
P_U(|u|)\le P_U(|\Pi_V(u)|)+P_U(P_V(|u|))=P_U(|\Pi_V(u)|)+P_V(|u|).
\]
By the definition of $\mathcal S$, $P_V(|u|)<\infty$ q.e. By \cite[Theorem 4.3.2]{FOT}, $P_U(|\Pi_V(u)|)\in F$, so it is finite q.e.
Consequently, $P_U(|u|)<\infty$ q.e. We also have  $\Pi_U(u)=\Pi_U(\Pi_Vu)$, which implies that $\Pi_U(u)\in F$.
Finally, by  properties of the orthogonal projection, for any $w\in F$ and $\eta\in F(U)$
we have $\EE(w,\eta)=\EE(\Pi_U(w),\eta)$. Therefore from \eqref{eq6.2} and Lemma \ref{lem2.5} we deduce that
\begin{equation}
\label{eq6.2ab}
\EE(\Pi_U(\Pi_V(u)),\eta)=\EE(\Pi_V(u),\eta)= \langle \mathbf1_{U}f(\cdot,u)\cdot m
+\mathbf1_{U}\cdot \mu,\eta\rangle,\quad \eta\in F(U).
\end{equation}
\end{remark}

\begin{definition}[Probabilistic solutions]
Let $\mu\in\RRr(D)$ and  $P_D|g|<\infty$.
We say that $u\in\mathcal B(E)$ is a {\em probabilistic solution} of (\ref{eq1.1}) if  $f(\cdot,u)\cdot m\in\RRr(D)$ and
\begin{equation}
\label{eq.probdef}
u=P_D (g)+R^Df(\cdot,u)+R^D\mu\quad \mbox{q.e.}
\end{equation}
\end{definition}

\begin{remark}
\label{rem.dqe}
Observe that by \eqref{eq2.26} and \eqref{eq2.8reg} we have that  \eqref{eq.probdef} holds if and only if $u=P_D (g)+R^Df(\cdot,u)+R^D\mu$ q.e. in $D$
and $u=g$ q.e. in $D^c$.
\end{remark}

\begin{theorem}
\label{th6.7}
Assume that $P_D(|g|)<\infty$ q.e.  and  $\mu\in\RRr(D)$.
\begin{enumerate}
\item[\rm(i)]Suppose that $u$ satisfies conditions (a),(b) of Definition \ref{def6.1} and
also (c) of this definition but only for $(V_n)$ appearing in (a).
Then  $u=P_D (g)+R^Df(\cdot,u)+R^D\mu$ q.e. in $D$.

\item[\rm(ii)] If $w\in\mathcal B(E)$, $f(\cdot,w)\cdot m\in \mathscr R(D)$ and $w=P_D (g)+R^Df(\cdot,w)+R^D\mu$ q.e., then $w$ is a  solution of  \eqref{eq1.1}.

\item[\rm(iii)] For any  solution $u$ of \eqref{eq1.1} and $V\in\mathcal O_q$,
if $\mathbf1_V|f(\cdot,u)|,\mathbf1_V\cdot|\mu|\in F^*$, then $\Pi_V(u)\in F$ and \eqref{eq6.2} holds with $V_n$ replaced by $V$.
\end{enumerate}
\end{theorem}
\begin{proof}
(i) Let   $u$ and $(V_n)$ be  as in (i).
Then, by condition (a) of Definition \ref{def6.1},
\[
\Pi_{V_n}(u)=R^{V_n}f(\cdot,u)+R^{V_n}\mu\quad\text{q.e. in }V_n.
\]
By this and conditions (b), (c) of Definition \ref{def6.1},
$u=P_D (g)+R^{D}f(\cdot,u)+R^D\mu$ q.e. in $D$. This proves (i).
Let  $w$ be as in (ii).
By  Lemma \ref{lem2.4}, (\ref{eq2.19})  and Corollary \ref{cor2.6}, for  every
$V\in \mathcal O_q$ such that $V\subset D$ we have
\begin{align*}
\Pi_V(w)&=\Pi_V(P_D(g))+\Pi_V(R^D\mu)+\Pi_V(R^{D}f(\cdot,w))\\
&=P_Dg-P_V(P_D((g))+R^{V}f(\cdot,u)+R^V\mu=R^{V}f(\cdot,w)+R^V\mu\quad \mbox{q.e.}
\end{align*}
By Lemma \ref{lem2.1} there is a sequence $(V_n)\subset\OO_q$ such that
$V_n\uparrow D$ q.e. and $\fch_{V_n}\cdot|\mu|, \fch_{V_n}f(\cdot,w)\cdot m\in F^*$, $n\ge1$. Set ${\mathcal S}=(V_n)$.
Then for every $V\in{\mathcal S}$, since $\mathbf1_V f(\cdot,w)$ and $\mathbf1_V\cdot\mu\in F^*$, we have
\begin{equation}
\label{eq6.8}
\Pi_V(w)=R^{V}f(\cdot,w)+R^V\mu=R^{V}(\mathbf1_Vf(\cdot,w))+R^V(\mathbf1_V\cdot\mu)\in F(V).
\end{equation}
Moreover, by (\ref{eq6.8}) and (\ref{eq2.20}), for every $V\in\mathcal S$,
\[
\EE(\Pi_V(w),\eta)=\EE(R^V(\mathbf1_V f(\cdot,w))+
R^V(\fch_V\cdot\mu),\eta)=\langle \mathbf1_V f(\cdot,w)+ \mathbf1_V\cdot\mu,\eta\rangle,\quad \eta\in F(V).
\]
Let $(U_n)$ be a $D$-total family satisfying (a) of Definition \ref{def6.1}.
By Corollary \ref{cor2.6} and (\ref{eq2.19}), for every $n\ge1$,
\begin{align*}
P_{U_n}(w)&=P_{U_n}(P_D(g))+P_{U_n}(R^Df(\cdot,w)+R^D\mu)\\&
=P_{D}(g)+(R^Df(\cdot,w)-\Pi_{U_n}(R^Df(\cdot,w)))+(R^D\mu-\Pi_{U_n}(R^D\mu))\quad \mbox{q.e.},
\end{align*}
so by Lemma \ref{lem2.4} and (\ref{eq2.12}),
\begin{align*}
P_{U_n}(w)-P_{D}(g)&=R^Df(\cdot,w)-R^{U_n}f(\cdot,w)+R^D\mu-R^{U_n}\mu\\&
=\mathbb E_\cdot \int_{\tau_{U_n}}^{\tau_D}f(X_s,w(X_s))\,ds+\mathbb E_\cdot (A^\mu_{\tau_D}-A^\mu_{\tau_{U_n}})\quad\text{q.e.}
\end{align*}
Since $\tau_{U_n}\nearrow \tau_D$ $\BP_x$-a.s. for q.e. $x\in D$, we see that
$P_{U_n}(w)-P_{D}(g)\to 0$ q.e. in $D$. Thus  $w$ is a projective variational solution of (\ref{eq1.1}),
which completes the proof of (ii).
Assertion (iii) is a consequence of  (i), (ii) and \eqref{eq6.8}.
\end{proof}

\begin{remark}
\label{rem.equiv}
The above proposition implies that in condition (c) of Definition \ref{def6.1}
one can  replace the word ``any''  by  ``some'' (cf. Remark \ref{rem.dqe}).
\end{remark}

\begin{remark}
Suppose that there exists a Green function for $L$ and $D$, i.e. a nonnegative Borel  function $G_D:E\times E\rightarrow\BR$ such that $G_D(x,y)=G_D(y,x)$, $G_D(x,y)=0$ if $x$ or $y$ belongs to $D^c$, $G_D(x,\cdot)$, $G_D(\cdot,y)$ are $(P^D_t)$-excessive for $x,y\in D$,  and moreover,
$R^Df(x)=\int_Df(y)G_D(x,y)\,m(dy)$ for any $x\in D$ and  bounded $f\in\BBr(D)$. Then for any $\mu\in\SSr(D)$ we have
\begin{equation}
\label{eq3.6}
R^D\mu(x)=\int_DG_D(x,y)\,\mu(dy)\quad \mbox{q.e. }x\in D.
\end{equation}
To show this one can argue as in the proof of  \cite[Lemma 3.1]{KR:CM} (applied to the form $\EE^D$). From (\ref{eq3.6}) and Theorem \ref{th6.7} it follows that $u$ is a probabilistic solution of (\ref{eq1.1}) if and only if it is an integral solution in the sense that (\ref{eq1.5green}) is satisfied for q.e. $x\in D$.
\end{remark}

\subsection{Existence and uniquenes of solutions.}

\begin{proposition}[Uniqueness result I]
\label{prop.unq1}
Assume \mbox{\rm(H1)}. There exists at most one  solution  $u\in L^1(D;m)$ of  \eqref{eq1.1}
such that $f(\cdot,u)\in L^1(D;m)$.
\end{proposition}
\begin{proof}
Let $u_1,u_2\in L^1(D;m)$ be  solutions of \eqref{eq1.1} such that we have $f(\cdot,u_1), f(\cdot,u_2)\in L^1(D;m)$.
Set $u=u_1-u_2$.  Let $\mathcal{S}_1=(V_n), \mathcal{S}_2=(W_n)$ be $D$-total families
of Definition \ref{def6.1} for $u_1$ and $u_2$, respectively. Write $T_1(w)(x)=\max\{\min\{w(x),1\},-1\}$ for any $w:E\rightarrow\BR$ and $x\in E$.
By \eqref{eq6.2} and Remark \ref{rem3.ext},
\begin{equation}
\label{eq6.unq1}
\EE(\Pi_{W_n\cap V_n}(u),T_1(\Pi_{W_n\cap V_n}(u)))=\int_D(f(\cdot,u_1)-f(\cdot,u_2))T_1(\Pi_{W_n\cap V_n}(u))\,dm.
\end{equation}
By \cite[Theorem 4.4.4]{FOT}, the form $\EE^D$ is transient. Hence, by \cite[Theorem 1.5.3]{FOT}, there exists a strictly positive function $\rho\in\mathcal B_b(D)$
and $c>0$ (depending only on $\rho$) such that
\[
\Big(\int_D|u|\rho\,dm\Big)^2\le c\EE^D(u,u),\quad u\in F(D).
\]
Moreover, for any Dirichlet form we have $\EE(u,T_1(u))\ge \EE(T_1(u),T_1(u))$, $u\in F(D)$.
Consequently,
\begin{equation}
\label{eq6.unq2}
\Big(\int_D|T_1(\Pi_{W_n\cap V_n}(u))|\rho\,dm\Big)^2\le c\int_D(f(\cdot,u_1)-f(\cdot,u_2))T_1(\Pi_{W_n\cap V_n}(u))\,dm.
\end{equation}
By \eqref{eq6.unq2}, condition (c) of Definition \ref{def6.1} and (H1) we get
\begin{equation}
\label{eq6.unq3}
\Big(\int_D|T_1(u)|\rho\,dm\Big)^2\le c\int_D(f(\cdot,u_1)-f(\cdot,u_2))T_1(u)\,dm\le 0.
\end{equation}
This proves the proposition.
\end{proof}

Applying the representation result of Theorem \ref{th6.7} and some probabilistic tools
we get a stronger than in Proposition \ref{prop.unq1} uniqueness result for solutions of \eqref{eq1.1}.
It  follows from the following comparison result.

\begin{theorem}
\label{prop6.6}
Let $f_1,f_2:E\times \mathbb R\to \mathbb R$ be Borel measurable functions such that $f_1(x,\cdot)$, $f_2(x,\cdot)$ are continuous for $x\in D$. Let
$g_1,g_2\in\mathcal B(E)$ be such that $P_D(|g_1|+|g_2|)<\infty$ q.e.,  and  let $\mu_1,\mu_2\in\RRr(D)$.
Assume that $\mu_1\le\mu_2$, $g_1\le g_2$  on $\partial_{\chi}D$,
$u_1,u_2$ are solutions of  \eqref{eq1.1} with $g,f,\mu$ replaced by $g_1,f_1,\mu_1$ and $g_2,f_2,\mu_2$, respectively,
and either $f_1$ is nonincreasing with respect to the second variable and $f_1(\cdot,u_2)\le f_2(\cdot,u_2)$ $m$-a.e. or
$f_2$ is nonincreasing with respect to the second variable and $f_1(\cdot,u_1)\le f_2(\cdot,u_1)$ $m$-a.e.
Then $u_1\le u_2$ q.e. in $D$.
\end{theorem}
\begin{proof}
By Theorem \ref{th6.7},
\[
u_i=P_D (g_i)+R^Df_i(\cdot,u_i)+R^D\mu_i\quad \text{q.e. in }D,\quad i=1,2.
\]
Let $w_i=u_i-h_i$, $h_i=P_D (g_i)$ and $\hat f_i(x,y)=f_i(x,y+h_i(x))$. Observe that
\[
w_i=R^D\hat f_i(\cdot,w_i)+R^D\mu_i\quad \text{q.e. in } D,\quad i=1,2.
\]
By Lemma \ref{lem4.2} and \cite[Proposition 4.9]{KR:JFA},
$w_1\le w_2$ q.e. in $D$, so  $u_1\le u_2$ q.e. in $D$.
\end{proof}

\begin{corollary}[Uniqueness result II]
\label{cor6.7}
Assume that \mbox{\rm(H1), (H3), (H4)} are satisfied. 
Then there exists at most one solution  of \eqref{eq1.1}.
\end{corollary}

\begin{lemma}
\label{lem6.9}
Let $\mu$ be a nonnegative smooth measure such that $\text{q-}\esssup_{D}R^D\mu<\infty$ (cf. \eqref{eq.a1}).
Then for q.e. $x\in E$,
\[
\mathbb E_x(A^\mu_{\tau_D})^2\le 2\|R^D\mu\|^2_\infty.
\]
\end{lemma}
\begin{proof}
By the strong Markov property and additivity of $A^\mu$,
\[
R^D\mu(X_t)=\mathbb E_{X_t}A^{\mu}_{\tau_D}=\mathbb E_x(A^\mu_{\tau_D}-A^\mu_t|\FF_t),\quad t\le\tau_D.
\]
Hence, by Lemma \ref{lm.itoinc}, for q.e. $x\in D$ we have
\[
\mathbb E_x(A^\mu_{\tau_D})^2=2\mathbb E_x\int_0^{\tau_D}\mathbb E_x(A^\mu_{\tau_D}-A^\mu_t|\FF_t)\,dA^\mu_t.
\]
Therefore
\[
\mathbb E_x(A^\mu_{\tau_D})^2=2\mathbb E_x\int_0^{\tau_D}R^D\mu(X_t)\,dA^\mu_t\le 2\|R^D\mu\|_\infty \mathbb E_xA^\mu_{\tau_D}\le 2\|R^D\mu\|^2_\infty,
\]
which proves the lemma.
\end{proof}

Let us note that by \cite[Theorem A.2.6, Theorem 4.1.3]{FOT},
for any $V\in\OO_q $,
\begin{equation}
\label{eq6.7}
X_{\tau_V}\circ\theta_{\tau_V}=X_{\tau_V}\quad \BP_x\text{-a.s.}\quad\text{for q.e. }x\in V.
\end{equation}
Let $\TT$ denote the  set of all stopping times with respect to the filtration $(\FF_t)_{t\ge0}$.
Note  also that an $(\FF_t)$-adapted c\`adl\`ag process is a  martingale with respect to the measure
$\BP_x$ if and only if for any bounded stopping time $\tau$ we have
$\mathbb E_x|M_\tau|<\infty$ and $\mathbb E_xM_\tau=\mathbb E_xM_0$.
A martingale $M$ with respect to $\BP_x$
is uniformly integrable if the family $\{M_\tau,\, \tau\in\TT, \tau<\infty\}$ is uniformly integrable.
The uniform integrability implies in particular that the limit $M_\infty=\lim_{t\to\infty}M_t$ $\BP_x$-a.s. exists.

\begin{proposition}
\label{prop6.10}
Let $V\in\OO_q $, $g\in\mathcal B(E)$ be such that $P_V(|g|)<\infty$ q.e, and $\mu$
be a nonnegative smooth measure such that $\mu\in\RRr(V)$.
Suppose that $u=P_V(g)+R^V\mu$ q.e. Then there exists a process $M$ such
that $M_0=0$, $M$ is a uniformly integrable martingale under the measure $\BP_x$
for q.e. $x\in V$ and
 \[
 u(X_t)=g(X_{\tau_V})+A^\mu_{\tau_V}-A^\mu_t
 -(M_{\tau_V}-M_t),\quad t\le \tau_V,\quad \BP_x\text{-a.s.}
 \]
 for q.e. $x\in V$.
\end{proposition}
\begin{proof}
Set $w=u-P_V(g)$. By \cite[Remark 3.3]{KR:NoD}, there exists a process $N$,
with the same properties as $M$ appearing in the assertion of the proposition and
such that for q.e. $x\in V$,
\[
w(X_t)=A^\mu_{\tau_V}-A^\mu_t-(N_{\tau_V}-N_t),\quad t\le \tau_V,\quad \BP_x\text{-a.s.}
\]
Let $h=P_V(g)$. We shall show that $L_t= h(X_t)-h(X_0)$, $t\le\tau_V$, shares the same properties as $N$.
Let $\alpha\in\TT,\, \alpha\le \tau_V$, $A=\{\alpha<\tau_V\}$, $B=\{\alpha=\tau_V\}$.
By the strong Markov property,
\begin{align*}
h(X_\alpha)=\mathbb E_{X_\alpha}g(X_{\tau_V})&=\mathbb E_x(g(X_{\tau_V}\circ\theta_\alpha)|{\FF_\alpha})\\
&=\mathbb E_x(\mathbf{1}_Ag(X_{\tau_V}\circ\theta_\alpha)|{\FF_\alpha})
+\mathbb E_x(\mathbf{1}_Bg(X_{\tau_V}\circ\theta_\alpha)|{\FF_\alpha}).
\end{align*}
On the set $A$ we have $\tau_V\circ \theta_\alpha=\tau_V-\theta_\alpha$,
so $\mathbf{1}_Ag(X_{\tau_V}\circ\theta_\alpha)=\mathbf{1}_Ag(X_{\tau_V})$.
Also, by (\ref{eq6.7}), $\mathbf{1}_Bg(X_{\tau_V}\circ\theta_\alpha)=\mathbf{1}_Bg(X_{\tau_V})$.
Hence $h(X_\alpha)=\mathbb E_x(g(X_{\tau_V})|\FF_\alpha)\, \BP_x$-a.s. for q.e. $x\in V$ and any $\alpha\in\TT$
such that $\alpha\le\tau_V$. As a result, the process $L$ has the required properties. Putting $M=N+L$ proves the proposition.
\end{proof}

\begin{theorem}
\label{th6.11}
Assume that \mbox{\rm(H1)--(H4)} are satisfied.
Then there exists a unique solution $u$ of \eqref{eq1.1}.
\end{theorem}
\begin{proof}
By \cite[Proposition 2.4, Theorem 2.9]{K:SPA}, for q.e. $x\in D$ there exists a unique stochastic process $Y^x$  such that
\begin{equation}
\label{eq5.8}
Y^x_t=g(X_{\tau_D})+\int_t^{\tau_D}f(X_s,Y^x_s)\,ds+A^\mu_{\tau_D}
-A^\mu_t-(M^x_{\tau_D}-M^x_t),\quad t\le\tau_D,\quad \BP_x\text{-a.s.}
\end{equation}
for some uniformly integrable martingale (with respect to $\BP_x$) $M^x$.
Moreover, by \cite[Theorem 2.9]{K:SPA} again, $\BE_x\int^{\tau_D}_0|f(X_t,Y^x_t)|\,dt<\infty$ for q.e. $x\in D$.
In view of (\ref{eq5.8}) and (\ref{eq2.11}), (\ref{eq2.25}), to prove the theorem  it is enough to show that
there exists a function $u:E\to\BR$ such that
\begin{equation}
\label{eq.aim1}
Y^x_t=u(X_t),\quad t\le\tau_D,\quad \BP_x\text{-a.s.}
\end{equation}
for q.e. $x\in D$. Indeed, if (\ref{eq.aim1}) holds true, then taking $t=0$ in (\ref{eq5.8})
and integrating with  respect to $\BP_x$ we get
\begin{align*}
u(x)=\mathbb E_x u(X_0)=\mathbb E_xY^x_0&=\mathbb E_xg(X_{\tau_D})+
\mathbb E_x\int_0^{\tau_D}f(X_t,u(X_t))\,dt+\mathbb E_xA^\mu_{\tau_D}\\
&=P_D (g)(x)+R^Df(\cdot,u)(x)+R^D\mu(x)
\end{align*}
for q.e. $x\in D$, so  by Theorem \ref{th6.7}, $u$ is a solution of \eqref{eq1.1}.
The proof of \eqref{eq.aim1} will be  divided  into three steps.
\\
{\em Step 1.} Suppose that there exists a strictly positive $\varrho\in\BB(E)$ such that $R^D\varrho<\infty$ q.e. and
\[
|f(x,y)|\le\varrho(x),\quad x\in E,\,y\in\BR.
\]
Let $h=P_D (g)$. By \cite[Theorem 4.7]{KR:JFA} there exists a unique $w\in\BB(E)$ such that
\[
w=R^Df_h(\cdot,w)+R^D\mu\quad\mbox{q.e.},
\]
where $f_h(x,y)= f(x,y+h(x))$, $x\in E,y\in\BR$. Let $u=w+h$. By Theorem \ref{th6.7},
$u$ is a solution of \eqref{eq1.1}. From \cite[Theorem 4.7]{KR:JFA} applied to $w$ and
Proposition \ref{prop6.10} applied to $h$ it follows that  $u(X)$ solves
\eqref{eq5.8}. By uniqueness, \eqref{eq.aim1} holds true.
\\
{\em Step 2.} Suppose now that
$\mathbb E_x(g(X_{\tau_D}))^2+\mathbb E_x(A^\mu_{\tau_D})^2+\mathbb E_x(\int_0^{\tau_D}|f(X_t,0)|\,dt)^2<\infty$.
Let $\varrho\in\BB(E)$ be a  strictly positive function such that $R^D\varrho<\infty$ q.e.
Such a function exists by \cite[Theorem 1.3.4.]{O}.
Set $\varrho_n= \frac{n\varrho}{1+n\varrho}$ and then
$f_{n,m}=\max\{\min\{f,n\varrho_n\},-m\varrho_m\}$ and $f_{m}=\max\{f,-m\varrho_m\}$.
By \cite[Lemma 2.7]{K:SPA},   for any $n,m\ge 1$ and q.e. $x\in D$ there exists a unique
$(\FF_t)$-adapted stochastic process  $Y^{x,n,m}$ (resp. $Y^{x,m}$) satisfying (\ref{eq5.8}),
with $f$ replaced by $f_{n,m}$ (resp. $f_m$) and $M^x$ replaced by a uniformly
integrable martingale $M^{x,n,m}$ (resp. $M^{x,m}$).
By the proof of \cite[Lemma 2.7]{K:SPA}, for q.e. $x\in D$ we have
\[
Y^{x,n,m}_t\nearrow Y^{x,m}_t,\quad t\le \tau_D,\quad \BP_x\mbox{-a.s.},
\]
and
\[
Y^{x,m}_t\searrow Y^x_t,\quad t\le \tau_D,\quad \BP_x\text{-a.s.}
\]
By {\em Step 1}, there exists a function $u_{n,m}$ on $E$ such that
$Y^{x,n,m}_t=u_{n,m}(X_t)$, $t\le\tau_D$, $\BP_x$-a.s. for q.e. $x\in D$.
By Theorem \ref{prop6.6}, $u_{n,m}\ge u_{n,m+1}, u_{n,m}\le u_{n+1,m}$ q.e. for all $n,m\ge 1$.
Set $u_m=\sup_{n\ge 1}u_{n,m}$ q.e. and $u=\inf_{m\ge 1} u_m$ q.e.
Then, by  \cite[Theorem 4.1.1, Theorem 4.2.1]{FOT},
\[
u_{n,m}(X_t)\nearrow u_m(X_t),\quad t\le \tau_D,\quad \BP_x\mbox{-a.s.},
\]
and
\[
u_{m}(X_t)\searrow u(X_t),\quad t\le \tau_D,\quad \BP_x\mbox{-a.s.}
\]
for q.e. $x\in D$. Combining the above convergences yields \eqref{eq.aim1}.
\\
{\em Step 3.} The general case. By \cite[Corollary 1.3.6]{O} applied to the form $\EE^D$ there exists a strictly positive
function $\varrho\in \BB(E)$ such that $R^D\varrho(x)<\infty$, $x\in D$.
Let $(V_k)\subset \mathcal O_q $ be a sequence such that $V_k\uparrow D$ q.e.
and $\|R^D (\mathbf1_{V_k}\cdot|\mu|)\|_\infty<\infty$, $k\ge 1$ (see Lemma \ref{lem2.1}).
Let $\mu_k=\mathbf1_{V_k}\cdot\mu$, $g_k=(g\wedge k)\vee (-k)$ and
\[
f_k(x,y)=f(x,y)-f(x,0)+(f(x,0)\wedge k)\vee(-k))\frac{k\varrho}{1+k\varrho}.
\]
By Lemma \ref{lem6.9}, the data $g_k, f_k, \mu_k$ satisfy the
assumptions of {\em Step 2}.  By \cite[Proposition 2.8]{K:SPA}, for any $k\ge 1$ and
q.e. $x\in D$, there exists a unique process $Y^{k,x}$  satisfying \eqref{eq5.8}
with $f$ replaced by $f_{k}$  and $M^x$ replaced by a uniformly martingale  $M^{x,k}$.
By the proof of  \cite[Proposition 2.8]{K:SPA}, for q.e. $x\in D$,
\begin{equation}
\label{eq.usf5}
\lim_{k\rightarrow\infty}\mathbb E_x\sup_{t\le\tau_D}|Y^{x,k}_t-Y^{x}_t|^{1/2}=0.
\end{equation}
On the other hand, by {\em Step 2}, there exists a function $u_{k}$ on $E$ such that
$Y^{x,k}_t=u_{k}(X_t)$, $t\le\tau_D$, $\BP_x$-a.s. for q.e. $x\in D$. It follows in
particular that $(u_k)$ is convergent q.e. in $D$.
Let $u=\lim_{k\to\infty}u_k$ q.e.
By  \cite[Theorem 4.1.1, Theorem 4.2.1]{FOT}, for q.e. $x\in D$ we have
\[
u_k(X_t)\to u(X_t),\quad t\le\tau_D,\quad \BP_x\mbox{-a.s.}
\]
as $k\rightarrow\infty$, which when combined with \eqref{eq.usf5} gives \eqref{eq.aim1}.
\end{proof}

We close this section with  a simple but important  corollary to Theorem \ref{th6.11}.

\begin{theorem}
\label{th3.hth}
Assume that \mbox{\rm(H1)--(H4)} are satisfied. Let
$h\in\mathcal B(E)$ be quasi integrable on $D$ and such that $R^D|f(\cdot,h)|<\infty$ q.e.
Then there exists a unique $u\in\mathcal B(E)$
such that $R^D|f(\cdot,u)|<\infty$ q.e. and
\begin{equation}
\label{eq1.5eh}
u(x)=h(x)+\mathbb E_x g(X_{\tau_D})+\mathbb E_x\int_0^{\tau_D}f(X_t,u(X_t))\,dt
+\mathbb E_x A^\mu_{\tau_D}\quad\text{q.e. }x\in E.
\end{equation}
\end{theorem}
\begin{proof}
For the existence part it is enough to apply Theorem \ref{th6.11} with $f$ replaced by
\[
f_h(x,y):=f(x,h(x)+y),\quad x\in E,\, y\in\mathbb R.
\]
Let $u_1,u_2$ be solutions to \eqref{eq1.5eh}. Observe that $u:=u_1-u_2$
is a solution of \eqref{eq1.1} with $\mu=0$, $g=0$
and $f=F$, where
\[
F(x,y):=f(x,y+u_2(x))-f(x,u_2(x)),\quad x\in E,\, y\in\mathbb R.
\]
Hence $u=0$ q.e. by Corollary \ref{cor6.7}.
\end{proof}

\subsection{Boundary trace operator.} Let us recall that by the  Beurling--Deny decomposition, for any  $u,v\in F$,
\begin{equation}
\label{eq1.bd}
\EE(u,v)=\EE^{(c)}(u,v)+\int_{E\times E}(u(x)-u(y))(v(x)-v(y))\,J(dx,dy)+\int_Eu(x)v(x)\kappa(dx),
\end{equation}
where $\EE^{(c)}$ is a symmetric form having the  strong local property,
$J$ is a symmetric Radon measure on $(E\times E)\setminus \mathfrak d$, where
$\mathfrak d:=\{(x,y)\in E\times E: x=y\}$, and $k$ is a positive Radon measure on $E$.
The above decomposition is unique (see, e.g., \cite[Lemma 4.5.4]{FOT}).
By $J^D,\kappa_D$ we denote the counterparts of $J,\kappa$ for the Dirichlet form $\EE$ restricted to $D$. Observe that
\[
\kappa_D(dx)=\mathbf1_D \cdot J(dx,D^c)+\mathbf1_D\cdot \kappa(dx).
\]

Now we are ready to formulate the second main result of this section.

\begin{theorem}
\label{th3.16}
\begin{enumerate}[\rm(i)]
\item Let $u$ be a solution of \eqref{eq1.1}. Then
for q.e. $x\in D$,
\begin{equation}
\label{eq1.wdnin}
\hat W^x_D(|u|):= \lim_{V\nearrow D, V\subset\subset D} P_V(|u|R^D\kappa_D)(x)=0.
\end{equation}

\item Let $u\in\mathcal B(E)$ be  bounded and
quasi continuous on $D\cup (\partial D\cap\partial_{\chi}D)$ and have the following properties:
\begin{enumerate}[\rm(1)]
\item$u$ satisfies conditions \mbox{\rm(a), (b)} of Definition \ref{def6.1},
\item $\hat W^x_D(u)=0$ q.e. in $D$.
\end{enumerate}
Then $u$ satisfies condition \mbox{\rm(c)} of Definition \ref{def6.1}.

\item Assume \mbox{\rm(H1), (H3), (H4)}. Then there exists at most one function $u\in\mathcal B(E)$ that is bounded,
quasi continuous on $D\cup (\partial D\cap\partial_{\chi}D)$ and  satisfies  conditions \mbox{\rm(1), (2)} of part \mbox{\rm(ii)}.
\end{enumerate}
\end{theorem}
\begin{proof}
In the proof we shall use the process $\mathbb M^D$ introduced before Lemma \ref{lem2.1}.

(i) By \cite[Lemma 4.5.2(iii)]{FOT} upon an application of the monotone convergence theorem and monotone class theorem; see \cite[Theorem I.8]{P}) and \cite[Theorem 4.4.2]{FOT},
\[
R^D\kappa_D(x)=\mathbb E^D_x1(X_{\tau_D-})=\mathbb E_x\mathbf1_D(X_{\tau_D-}),\quad x\in D.
\]
Hence, by the strong Markov property,
\[
 \mathbf1_{\{\tau_V<\tau_D\}}R^D\kappa_D(X_{\tau_V})= \mathbf1_{\{\tau_V<\tau_D\}}\mathbb E^D_{X_{\tau_V}}1(X_{\tau_{D}-})
 =\mathbb E_x^D\big( \mathbf1_{\{\tau_V<\tau_D\}} 1(X_{\tau_D-})\big|\FF_{\tau_V}\big)
 \quad \BP^D_x\mbox{-a.s.}
\]
for q.e. $x\in D$. Consequently,
\begin{align}
\label{eq.spl1}
P_V(|u|R^D\kappa_D)(x)&=\mathbb E_x(|u|(X_{\tau_V})R^D\kappa_D(X_{\tau_V}))\nonumber\\
&=\mathbb E_x\big([\mathbf1_D|u|](X_{\tau_V})\cdot[R^D\kappa_D](X_{\tau_V})\big)\nonumber\\
&=\mathbb E^D_x\big(\mathbf1_{\{\tau_V<\tau_D\}}
[\mathbf1_D|u|](X_{\tau_V})\cdot[R^D\kappa_D](X_{\tau_V})\big)
\nonumber\\&=
\mathbb E^D_x\big(\mathbf1_{\{\tau_V<\tau_D\}}[\mathbf1_D|u|](X_{\tau_V})\cdot 1(X_{\tau_D-})\big)\nonumber \\
&= \mathbb E^D_x[|u|(X_{\tau_V})\mathbf1_{\{X_{\tau_D-}\in D\}}]
=\mathbb E_x[\mathbf1_D|u|(X_{\tau_V})\mathbf1_{\{X_{\tau_D-}\in D\}}].
\end{align}
Let $(V_n)$ be an increasing sequence or relatively compact open subsets of $D$
such that $V_n\uparrow D$ q.e.  Since $X$ is quasi-left continuous under $\BP_x$, we have
\begin{equation}
\label{eq.fmt1}
\BP_x(\bigcup_{n\ge 1}\{\tau_{V_n}=\tau_D\}\cap\{X_{\tau_D-}\in D\} ) =\BP_x(\{X_{\tau_D-}\in D\}) \quad \text{q.e. }x\in D.
\end{equation}
As a result,
\[
\mathbf1_D|u|(X_{\tau_{V_n}})\mathbf1_{\{X_{\tau_D-}\in D\}}\to 0\quad \BP_x\text{-a.s.}\,\, \text{for q.e. }x\in D.
\]
What is left is to  show that $(|u|(X_{\tau_{V_n}}))$ is uniformly integrable under
the measure $\BP_x$
for q.e. $x\in D$. But this follows from Lemma \ref{lem4.2}. Thus \eqref{eq1.wdnin} holds.

(ii) Let $u\in\mathcal B(E)$ be   bounded and quasi continuous  on $D\cup (\partial D\cap\partial_{\chi}D)$, and
 satisfy (1), (2).
We shall prove that (c) of Definition \ref{def6.1} holds true.
Let $(V_n)$ be a $D$-total family such
that $V_n\subset\subset D$, $n\ge 1$ (see Theorem \ref{th6.7}(i)). We have
\begin{align*}
P_{V_n}u(x)&=\mathbb E_x\big[(\mathbf1_{D}u)(X_{\tau_{V_n}})
\mathbf1_{\{X_{\tau_D-}\in D\}}\big]+
\mathbb E_x\big[(\mathbf1_{D^c}u)(X_{\tau_{V_n}})\mathbf1_{\{X_{\tau_D-}\in D\}}\big]\\
&\quad+
\mathbb E_x\big[u(X_{\tau_{V_n}})\mathbf1_{\{X_{\tau_D-}\notin D\}}\big]\\
&=P_{V_n}(uR^D\kappa_D)(x)
+\mathbb E_x\big[(\mathbf1_{D^c}g)(X_{\tau_{D}})\mathbf1_{\{X_{\tau_D-}\in D\}}\mathbf1_{\{\tau_{V_n}=\tau_D\}}\big]\\&
\quad+
\mathbb E_x\big[u(X_{\tau_{V_n}})\mathbf1_{\{X_{\tau_D-}\notin D\}}\big],
\end{align*}
where the second equality being a consequence of  \eqref{eq.spl1}. By the assumptions we made, we have
$P_{V_n}(uR^D\kappa_D)(x)\to \hat W^x_D(u)=0$
q.e. in $D$. Observe also that
\[
\mathbb E_x\big[(\mathbf1_{D^c}g)(X_{\tau_{D}})\mathbf1_{\{X_{\tau_D-}\in D\}}\mathbf1_{\{\tau_{V_n}=\tau_D\}}\big]
\to \mathbb E_x\big[(\mathbf1_{D^c}g)(X_{\tau_{D}})\mathbf1_{\{X_{\tau_D-}\in D\}}\big]
\]
(see \eqref{eq.fmt1}). Now, if $X_{\tau_D-}\notin D$ then $\tau_{V_n}<\tau_D,\, n\ge 1$ and
$\lim_{n\to \infty} X_{\tau_{V_n}}=X_{\tau_D}\in \partial D\cap \partial_{\chi}D$ (we use
quasi-left continuity of $\BM$,  which implies that $X_{\tau_{V_n}}\to X_{\tau_D}$ $\BP_x$ a.s. for q.e. $x\in D$).
Since $u$ was assumed to be quasi continuous on $D\cup (\partial D\cap \partial_{\chi}D)$, applying \cite[Theorem 4.2.2]{FOT} shows that for q.e. $x\in D$ we have
 $u(X_{\tau_{V_n}})\to u(X_{\tau_D})=g(X_{\tau_D})$ $\BP_x$-a.s. on the set $\{X_{\tau_D-}\notin D\}$.
This when combined with the assumption that  $u$ is bounded on $D\cup (\partial D\cap \partial_{\chi}D)$
implies, by the Lebesgue dominated convergence theorem, that
\[
\mathbb E_x\big[u(X_{\tau_{V_n}})\mathbf1_{\{X_{\tau_D-}\notin D\}}\big]\to
\mathbb E_x\big[g(X_{\tau_{D}})\mathbf1_{\{X_{\tau_D-}\notin D\}}\big]
\]
for q.e. $x\in D$. Putting all the convergences together we see that $P_{V_n}u\to P_Dg$ q.e. in $D$.
This finishes the proof of (ii). Assertion (iii) is a consequence of (ii) and Corollary \ref{cor6.7}.
\end{proof}

\section{Purely nonlocal operators and solutions with nonzero boundary trace}
\label{sec7ab}

In the this section, we  focus on semilinear equations with $L$ belonging to
a   special class of purely nonlocal operators.
Throughout this section,
$E=\mathbb R^d$ and  the Dirichlet form $\EE$ is assumed to bo purely jumping, i.e.
$\EE^{(c)}=0$ in the Beurling--Deny decomposition \eqref{eq1.bd}.
We also assume that $\kappa=0$.

\begin{lemma}
\label{lm.harq}
Assume that $\BP_x(\tau_D<\infty)=1$ q.e. and $h$ is an $m$-a.e. finite $(P^D_t)$-excessive function.
Then $h$ is quasi integrable on $D$.
\end{lemma}
\begin{proof}
Since $h\in\mathbb W(D)$, we have $R^D_1h(x)\le h(x),\, x\in D$.
Since   $h\in\mathbb W(D)$ and is finite $m$-a.e. in $D$, it is finite q.e. in $D$ (see \cite[Theorem A.2.13]{CF}).
Thus
\[
\mathbb E_x\int_0^{\tau_D}e^{-t}h(X_t)\,dt<\infty\quad \text{q.e. }x\in E.
\]
From this we readily get the result.
\end{proof}

Throughout the rest of the section we assume that $D$ is bounded and  the following condition holds.
\begin{enumerate}
\item[(B)] $L$ is  of the form
\[
L=-\phi(-\Delta),
\]
where $\phi:(0,\infty)\to [0,\infty)$ is  a complete Bernstein function
with  L\'evy density  $\hat\mu:[0,\infty)\to [0,\infty)$, i.e.
\[
\phi(\lambda)=\int_0^\infty (1-e^{-\lambda t})\hat \mu(t)\,dt,\quad \lambda>0,
\]
(by the definition, $\int_0^\infty(1\wedge t)\hat\mu(t)\,dt<\infty$).
Moreover, there exist $a_1,a_2, R_0>0$ and $0<\delta_1\le\delta_2<1$ such that
\[
a_1\Big(\frac{t}{s}\Big)^{\delta_2}\le \frac{\phi(t)}{\phi(s)}\le a_2\Big(\frac{t}{s}\Big)^{\delta_2},\quad R_0\le s\le t.
\]
\end{enumerate}

By \cite{Bi} there exists Green's function $G_D$ for $L_D$
and for any $y\in \partial D$ the limit of $M_D(x,y):=G_D(x,y)/G_D(x_0,y)$
exists as $D\ni y\to y_0$. Therefore $M_D(x,y)$ (the so-called Martin kernel) is well defined for $x\in D,\, y\in\bar D$.
It is known that under condition  (B)  the operator $L$  admits  the form \eqref{eq1.11} with
\[
j(r)=\frac{1}{(4\pi t)^{d/2}}\int_0^\infty e^{-r^2/4t}\hat\mu(t)\,dt.
\]
Let $\partial_mD:=\partial D\setminus \{x\in\partial D: x\text{ is inaccessible from }D\}$.
Recall that a point $x\in\partial D$ is called inaccessible from $D$ if for any $x_0 \in D$,
\[
\int_D G_D\left(x_0, z\right) j(|z-y|) d z<\infty.
\]
For any nonnegative Borel measure $\nu$ on $\partial_m D$ we let
\[
M_D\nu(x)=\int_{\partial_mD}M_D(x,y)\,\nu(dy),\quad x\in D.
\]
By \cite[Theorem 4.3]{Bi}, $M_D\nu$ is harmonic in $D$.
We start with showing that \eqref{eq.ident1} holds.
Fix $x_0\in D$.
\begin{proposition}
For any $A\in\mathcal B(\mathbb R^d)$ and any
Lipschitz regular open set $V\subset\subset D$ we have
\[
\eta_V[u](A)=\mathbb E_{x_0}(\mathbf1_Du(X_{\tau_V})\mathbf1_A(X_{\tau_V-})),
\]
where $\eta_V[u]$ is defined by \eqref{eq1.12}. In particular, \eqref{eq.ident1} holds true.
\end{proposition}
\begin{proof}
Observe that
\begin{align*}
\mathbb E_{x}(\mathbf1_Du(X_{\tau_V})\mathbf1_A(X_{\tau_V-}))
&=
\mathbb E_{x}(\mathbf1_{D\setminus V}u(X_{\tau_V})\mathbf1_{A\cap V}(X_{\tau_V-}))
\\&=
\mathbb E_{x}(\mathbf1_{D\setminus \bar V}u(X_{\tau_V})\mathbf1_{A\cap V}(X_{\tau_V-})),
\end{align*}
where in the last equation we used regularity of $V$.
By \cite[Lemma 4.5.5]{FOT}, for any $h\in\mathcal B_b^+(V)$,
\begin{align*}
\mathbb E_{h}(\mathbf1_{D\setminus \bar V}u(X_{\tau_V})\mathbf1_{A\cap V}(X_{\tau_V-}))
&= 2\int_{\mathbb R^d\times\mathbb R^d} R^Vh(x)
\mathbf1_{A\cap V}(x)\mathbf1_{D\setminus \bar V}u(y)\,J(dx,dy)\\
& =\int_{A\cap V}\int_{D\setminus \bar V} R^Vh(x)u(y)j(|x-y|)\,dy\,dx
\\&
= \int_{A\cap V}\int_{D\setminus \bar V} \int_V G_V(x,z)h(z)\,dz \,u(y)j(|x-y|)\,dy\,dx.
\end{align*}
Hence, for a.e. $z\in V$,
\begin{align*}
\mathbb E_{z}(\mathbf1_Du(X_{\tau_V})\mathbf1_A(X_{\tau_V-}))
= \int_{A\cap V}\int_{D\setminus \bar V} G_V(z,x) \,u(y)j(|x-y|)\,dy\,dx.
\end{align*}
Since both sides of the above equality are $(P^V_t)$-excessive,
it holds for every $z\in V$. This completes the proof of the first assertion.
By what has already been proved and (\ref{eq2.25}) we have  $\eta_V[u](\BR^d)=\mathbb E_{x_0}(\fch_Du(X_{t_V}))=P_V(\fch_Du)(x_0)$. To get \eqref{eq.ident1}
it suffices now to observe that $R^D\kappa_D=\mathbf1_D$ and use (\ref{eq1.12}) and \eqref{eq1.wdnin}.
\end{proof}

We let
\begin{equation}
\label{eq.pkac}
p_D(x,y):=\int_DG_D(x,z)j(|z-y|)\,dz,\quad x\in D,\, y\in D^c.
\end{equation}

\begin{theorem}
\label{th4.3}
Let $\nu$ be a bounded Borel measure on $\partial_mD$,  $\gamma$ be a
Borel measure on $\mathbb R^d\setminus(D\cup\partial_mD)$
and  $P_D(|\gamma|)<\infty$ q.e. in $D$.
Assume that \mbox{\rm(H1)--(H4)} are satisfied and
\[
\int_D |f(y,M_D\nu(y))|G_D(x,y)\,dy<\infty\quad \text{q.e. }x\in D.
\]
Then there exists a unique function $u$ such that
\begin{align*}
u(x)&=\int_{\partial_m D}M_D(x,y)\,\nu(dy)+\int_{(D\cup\partial_mD)^c}p_D(x,y)\,\gamma(dy)\\
&
\quad+\int_Df(y,u(y))G_D(x,y)\,dy+\int_DG_D(x,y)\,\mu(dy)\quad \text{q.e. }x\in D.
\end{align*}
Furthermore,
\[
W_D[u]=\nu.
\]
\end{theorem}
\begin{proof}
The existence part follows  from Lemma \ref{lm.harq} and Theorem \ref{th3.hth}
with  $h=M_D\nu$, $g=0$ and $\mu$ replaced by $\mu+\beta$, where
\[
\beta(x)= \int_{D^c}j(x,y)\gamma(dy).
\]
Note that by an easy application of Fubini's theorem we have
\[
R^D\beta(x)=\int_{(D\cup\partial_mD)^c}p_D(x,y)\,\gamma(dy)\quad\text{q.e. }x\in D.
\]
The second assertion is a consequence of \cite[Propositions  5.4, 5.11]{Bi}.
\end{proof}

\section{Regularity results I and Sobolev spaces in the broad sense}

In this section, we shall  introduce some spaces
which seem to be natural when studying regularity of \eqref{eq1.1}
with nonregular data.

It is well know that even for the classical Dirichlet problem, i.e.
when $D$ is a smooth bounded domain, $L=\Delta$, $f=\mu=0$ and $g\in C(\partial D)$,
in general, the solution $u$ of \eqref{eq1.1} is only locally in the energy space,
i.e. $u\in H^1_{loc}(D)$.
In addition, if $\mu$ is nontrivial, then in general $u\notin H^1_{loc}(D)$,
and the best regularity one can get  is that $T_k(u):=\max\{\min\{u,k\},-k\}\in H^1_{loc}(D)$
for any $k\ge 1$. This  means that $u\in H^1_{loc}(\{|u|<k\})$  for $k\ge 1$.
In general, $\{|u|<k\}$ is not open but only quasi open. Therefore  in our framework it is natural  to consider Sobolev spaces on quasi open sets.   However, the notion of  a``local property'' for functions  that solve equations with  nonlocal operators is not such a natural and straightforward concept.
Our goal in this section  is to give some definition of local energy spaces for (possibly) nonlocal operators.

For a  family $\mathcal S\subset \mathcal O_q(W)$, we define
\[
F_{\chi}(W;\mathcal S):=\{u\in\BB^n(E):  P_V(|u|)< \infty\,\,m\text{-a.e. and }\Pi_V(u)\in F(V),\, V\in \mathcal{S}\},
\]
and
\[
F(W;\mathcal S)=\{u\in\BB^n(E):  \text{for each }V\in \mathcal{S}\text{ there is } \eta\in F \text{ such that } u=\eta\text{ q.e. on } V\}.
\]

For $W\in\OO_q $ we define
\[
\Xi_W=\{\mathcal S:\mathcal S\mbox { is a $W$-total family}\},
\]
and then
\[
\dot{F}_{\chi,\text{loc}}(W)=\bigcup_{\mathcal S\in\Xi_W} F_{\chi}(W;\mathcal S),\qquad
\dot{F}_{\text{loc}}(W)=\bigcup_{\mathcal S\in\Xi_W} F(W;\mathcal S).
\]

\begin{remark}
\label{rem3.2}
Let $u\in F_{\chi}(W;\mathcal S)$ and $V\in \mathcal S$.
By the definition, $u-P_V(u)=\eta$ for some $\eta\in F(V)$.
It is a matter of straightforward calculation that $P_V(|u|)\in\mathbb W(D)$,
hence $P_V(|u|)$ is $\EE^V$-quasi continuous. Consequently $u=P_V(u)+\eta$
is $\EE^V$-quasi continuous. Since the family $\mathcal S$ is $W$-total,
we conclude that $u$ is $\EE^W$-quasi continuous.
\end{remark}

Note that  the space   $\dot F_{loc}(E)$ is considered
in \cite[p. 271]{FOT} and \cite[page 163]{CF}.
As in  \cite{FOT},  the elements of the space $\dot F_{loc}(W)$
may be called functions which are locally in $F(W)$ in the broad sense.
By introducing  the space $\dot{F}_{\chi,\text{loc}}(W)$ we want to express
in a different manner, when comparing to $\dot{F}_{\text{loc}}(W)$,
the fact that ``$u$ is in the energy space $F$ on  parts $V$ of $W$''.
Instead of demanding that $u$ may be extended from $V$ to $E$ in such a way that the extension belongs to $F$,
we demand that its projection $\Pi_V(u)$ belongs to $F(V)\subset F$.
This property better corresponds to the definition of a solution of \eqref{eq1.1}.

Let $V\in\mathcal O_q$. In what follows we denote by  $\mbox{Cap}_{\EE^V}$
the capacity associated with $\EE^V$.

\begin{lemma}
\label{lem5.2}
Let $V\in\OO_q $ and $\CC_V=\{U\in\OO_q : U\subset V,\, \mbox{\rm Cap}_{\EE^V}(U)<\infty\}$. Then
there exists an increasing sequence $(U_n)_{n\ge1}\subset \mathcal C_V$ such that
\[
\bigcup_{n\ge 1}U_n=V\quad\text{q.e.}
\]
\end{lemma}
\begin{proof}
By \cite[Corollary 1.3.6]{O} there exists a strictly positive quasi continuous
function $g\in F$
such that $Rg(x)\le 1$, $x\in E$.  Observe that $R^Vg\in F(V)$ and $R^Vg(x)>0$, $x\in V$.
Set $U_n=\{R^Vg>1/n\}$. Then clearly $\bigcup_{n\ge 1}U_n=V$ and $nR^Vg\ge \mathbf1_{U_n}$
q.e. in $V$. Since $nR^Vg\in F(V)$, we deduce that $\mbox{Cap}_{\EE^V}(U_n)<\infty$.
Since $(U_n)_{n\ge 1}\subset \mathcal C_V$, we get the result.
\end{proof}

\begin{theorem}
\label{th5.3}
Suppose that $V\in\mathcal O_q$  and $\Pi_V(u)\in F$.
Then for any $U\in\mathcal C_V$ such that $P_V(|u|)\le c$ q.e. in $U$
there exists $\eta_U\in F$ such that
$u=\eta_U$ q.e. in $U$. Consequently, for any $W\in \mathcal O_q$,
\[
\dot F_{\chi,\text{loc}}(W)\subset \dot F_{\text{loc}}(W).
\]
\end{theorem}
\begin{proof}
Let $V$, $u$  satisfy the assumptions of the theorem.
Then $\Pi_V(u)\in F(V)$, so there exists $\xi\in F(V)$ such that $u=P_V(u)+\xi$ q.e.
Let $U\in\OO_q$
be such that $U \subset V$,  $\mbox{Cap}_{\EE^V}(U)<\infty$ and $P_V(|u|)\le c$ q.e. in $U$.
\\
{\em Step 1.} We shall show that there exists $G\in\OO_q$ such that $U\subset G\subset V$
and $\mbox{Cap}_{\EE^G}(U)<\infty$, $\mbox{Cap}_{\EE^V}(G)<\infty$. Since $\mbox{Cap}_{\EE^V}(U)<\infty$,
there exists the 0-equilibrium potential $e^V_U\in F(V)$ (with respect to the form $\EE^V$) and
$\delta\in\SSr(V)$ (the 0-equilibrium measure) such that  $e^V_U=R^V\delta$ q.e. in $V$ (see \cite[p. 82]{FOT}
and the comments following \cite[Corollary 2.2.2]{FOT}).
Set $G=\{R^V\delta>1/2\}$. Since $R^V\delta$ is quasi continuous, $G$ is quasi open, and since
$R^V\delta=1$  on $U$, we have $U\subset G$.
Since every normal contraction operates on $\EE^V$,  $f:=2(R^V\delta-1/2)^+\in F(G)$.
Observe that $f\ge 1$ q.e. on $U$. Hence $\mbox{Cap}_{\EE^G}(U)<\infty$.
On the other hand, $2R^V\delta\ge 1$ q.e. on $G$. Hence  $\mbox{Cap}_{\EE^V}(G)<\infty$.
\\
{\em Step 2.}
Let $U,G$ be as in {\em Step 1}.
Let  $g= e^{G}_{U}$. Set $h=P_V(u)$, $h_1=P_V(u^+)$ and $h_2=P_V(u^-)$.
We then have
\begin{equation}
\label{xih}
u=\xi +h.
\end{equation}
Let $e^G_{U,h_1}, e^G_{U,h_2}\in F(G)$ be such that $0\le e^G_{U,h_i}\le h_i$
on $G$ and $e^G_{U,h_i}=h_i$ on $U$, $i=1,2$ (see \cite[Exercise 3.10]{MR}).
Observe that $e^G_{U,h_i}g=h_i$ q.e. on $U$, $i=1,2$,  so $(e^G_{U,h_1}-e^G_{U,h_2})g=h$ q.e. on $U$. Since $e^V_{G,h_i}, g\in F$ and each of them is bounded, $e^V_{G,h_i}g\in F$, $i=1,2$, by  \cite[Corollary 1.5.1]{FOT}. This finishes the proof of the first assertion. The second one is a consequence of the first one, Lemma \ref{lem5.2}
and the fact that $\{|u|<k\}\cap W$ is a $W$-total family.
\end{proof}


\begin{corollary}
Let   $u$ be a solution of  \mbox{\rm(\ref{eq1.1})}. Then    $u\in\dot F_{\chi,\text{loc}}(D)\subset \dot F_{\text{loc}}(D)$.
\end{corollary}

\section{Regularity results II and  a priori estimates}
\label{sec7}

In Section \ref{sub2.1} we have introduced the space $\MMr_\rho(D)$. We equip it with the norm
\[
\|\mu\|_{\MMr_\rho(D)}:=\int_D\rho\,d|\mu|<\infty.
\]

\begin{proposition}
\label{prop7.1}
Assume \mbox{\rm (H1)--(H3)}. Let $u$ be a solution of \eqref{eq1.1}. Then
\begin{enumerate}[\rm(i)]
\item For q.e. $x\in D$
\begin{equation}
\label{eq7.1}
|u|+R^D|f(\cdot,u)|\le  2R^D|f(\cdot,0)|+R^D|\mu|+P_D|g|.
\end{equation}

\item For q.e. $x\in D$,
\[
|u-P_Dg|+R^D|f(\cdot,u)|\le  2R^D|f(\cdot,P_Dg)|+R^D|\mu|,
\]

\item For any $\rho\in\mathbb W(D)$ we have
\[
\|f(\cdot,u)\|_{L^1_\rho(D;m)}\le  2\|f(\cdot,P_Dg)\|_{L^1_\rho(D;m)}+\|\mu\|_{\MMr_\rho(D)}.
\]
\end{enumerate}
\end{proposition}
\begin{proof}
(i) Let $u$ and $M^x$ be defined as in the proof of Theorem \ref{th6.11}, and let $Y=u(X)$.
By \cite[, Proposition 2.4, Theorem 2.9]{K:SPA}, $Y$ is of Doob's class (D) under $\BP_x$ for q.e. $x\in D$, i.e. the collection of random variables $\{Y_{\tau},\tau\in\TT,\tau<\infty\}$ is uniformly integrable under the measure $\BP_x$ for q.e. $x\in D$. Furthermore, by
(\ref{eq5.8}) and (\ref{eq.aim1}), for q.e. $x\in D$ we have
\begin{align}
\label{eq7.2}
Y_t&=g(X_{\tau_D})+\int^{\tau_D}_{t\wedge\tau_D}f(X_s,u(X_s))\,ds\nonumber \\
&\quad+A^{\mu}_{\tau_D}-A^{\mu}_{t\wedge\tau_D}-(M^x_{\tau_D}-M^x_{t\wedge\tau_D}), \quad t\ge0,\quad \BP_x\mbox{a.s.}
\end{align}
Applying the Meyer--Tanaka formula  (see \cite[Corollary 3 to Theorem IV.70]{P}) we get
\begin{align*}
|Y_{t\wedge\tau_D}|-|Y_0|&\ge\int^{t\wedge\tau_D}_0
\mbox{sgn}(Y_{s-})\,dY_s=-\int^{t\wedge\tau_D}_0\mbox{sgn}(Y_s)(f(X_s,Y_s)-f(X_s,0))\,ds\\
&\quad+\int^{t\wedge\tau_D}_0\mbox{sgn}(Y_s)(-f(X_s,0)\,ds+dA^{\mu}_s)
+\int^{t\wedge\tau_D}_0\mbox{sgn}(Y_{s-})\,dM^x_s,
\end{align*}
where $\mbox{sgn}(x)=1$ if $x>0$ and $\mbox{sgn}(x)=-1$ if $x\le0$. By this and (H1),
\begin{align}
\label{eq7.3}
&|u(x)|+\BE_x\int^{t\wedge\tau_D}_0|f(X_s,u(X_s))-f(X_s,0)|\,ds\nonumber\\
&\qquad\le \BE_x|u(X_{t\wedge\tau_D})|+\BE_x\int^{t\wedge\tau_D}_0(|f(X_s,0)|\,ds
+dA^{|\mu|}_s),\quad t>0.
\end{align}
By (\ref{eq7.2}), for q.e. $x\in D$, $u(X_{t\wedge\tau_D})\rightarrow g(X_{\tau_D})$ $\BP_x$-a.s. as $t\rightarrow\infty$.
Since $u(X)$ is of class (D), it follows that   $\BE_x|u(X_{t\wedge\tau_D})|\rightarrow \BE_x|g(X_{\tau_D})|$ as $t\rightarrow\infty$.
Therefore letting $t\rightarrow\infty$ in (\ref{eq7.3}) yields (\ref{eq7.1}).
In order to get (ii) it is enough to observe that $w=u-P_Dg$
is a solution of \eqref{eq1.1} with $g=0$ and $f$ replaced by $f_g(x,y):=f(x,y+P_Dg(x))$.
Applying (i) to $w$ yields (ii). (iii) is a consequence of (ii)
and \cite[Lemma 4.6]{K:CVPDE}.
\end{proof}

\begin{corollary}
Assume \mbox{\rm (H1)--(H3)}. Let $u$ be a solution of \eqref{eq1.1}.
\begin{enumerate}[\rm(i)]
\item If  $(W_n)\subset\mathcal O_q$ is $D$-total family  such that
\[
\int_{W_n}\max\{|f(\cdot,n)|,|f(\cdot,-n)|\}\,dm
+\int_{W_n}P_D|g|\,dm+|\mu|(W_n)<\infty,\quad n\ge 1,
\]
then $\Pi_{V_n}(u)\in F$, $n\ge 1$, where
$V_n=\{R^D|f(\cdot,0)|+R^D|\mu|+P_D|g|<n\}\cap W_n$.

\item If  $(U_n)\subset\mathcal O_q$ is a $D$-total family
and
\[
\int_{U_n}|f(\cdot,P_Dg)|\,dm+|\mu|(U_n)<\infty,\quad n\ge 1,
\]
then $\Pi_{\hat U_n}(u)\in F$, $n\ge 1$, where  $\hat U_n
=\{R^D|f(\cdot,P_Dg)|+R^D|\mu|<n\}\cap U_n$.
\end{enumerate}
\end{corollary}
\begin{proof}
By \eqref{eq.probdef} and Lemma \ref{lem2.4},   for any $V\in \mathcal O_q$ we have
\[
\Pi_V(u)=R^Vf(\cdot,u)+R^V\mu.
\]
By \cite[Proposition 5.9]{KR:JFA}, if $R^V|f(\cdot,u)|+R^V|\mu|$ is bounded
and $|\mu|(V)+\int_V|f(\cdot,u)|\,dm<\infty$, then $R^Vf(\cdot,u), R^V\mu\in F(V)$,
and hence $\Pi_V(u)\in F(V)\subset F$. From this,  the choice of $V_n$ and $\hat U_n$ and
Proposition \ref{prop7.1} we conclude the result.
\end{proof}

\begin{proposition}
\label{prop7.2}
Assume that $(f_1,g_1,\mu_1), (f_2,g_2,\mu_2)$ satisfy \mbox{\rm(H1)--(H3)}.
Let $u_i$, $i=1,2$, be a solution of \eqref{eq1.1}
with $(f,g,\mu)$ replaced by $(f_i,g_i,\mu_i)$.  Then
\[
|u_1-u_2|\le R^D|f_1(\cdot,u_1)-f_2(\cdot,u_1)|
+R^D|\mu_1-\mu_2|+P_D|g_1-g_2|
\]
for q.e. in $D$. If, in addition $f_1=f_2=:f$, then
\[
|u_1-u_2|+R^D|f(\cdot,u_1)-f(\cdot,u_2)|\le R^D|\mu_1-\mu_2|+P_D|g_1-g_2|.
\]
\end{proposition}
\begin{proof}
It suffices to observe that $u_1-u_2$ is a solution to \eqref{eq1.1}
with $g$ replaced by $g_1-g_2$, $f$ replaced by $F(x,y):=f_1(x,y+u_2(x))-f_2(x,u_2(x))$,
and $\mu$ replaced by $\mu_1-\mu_2$, and then apply Proposition \ref{prop7.1}.
\end{proof}

In the examples below, $D$ is a bounded open subset of  $\BR^d$, $d\ge3$. We set
\[
\delta(x)=\inf_{y\in\partial D}|x-y|,\quad x\in\BR^d.
\]
For nonnegative real functions $u,v$ on $D\subset E$ the notation $u\asymp v$
means that $c^{-1}u\le v\le c u$ for some $c\ge1$.

\begin{example}
\label{ex7.6}
Assume additionally that $D$ is of class $C^2$. Let $m$ (or $dx$) denote the Lebesgue measure
on $D$ and $\sigma$ denote the surface measure on $\partial D$. Consider equation (\ref{eq1.1})
with $L=\Delta$ (see Example \ref{ex2.7}(ii)). It is well known (see \cite{MV,Z})
that the  Green function $G_D$ and the density $p_D(x,y)$ of the  Poisson kernel satisfy the following estimates
\begin{equation}
\label{eq7.4}
G_D(x,y)\asymp \min\big\{|x-y|^{2-d}\,,
\delta(x)\delta(y)|x-y|^{-d}\big\},\quad x,y\in D,
\end{equation}
and
\begin{equation}
\label{eq7.5}
p_D(x,y)\asymp \delta(x)|x-y|^{-d},\quad x\in D,\,y\in\partial D.
\end{equation}
As in the proof of \cite[Proposition 4.9]{Ku} one can show that from (\ref{eq7.4}) it follows that
\begin{equation}
\label{eq7.6}
R^D1(x)=\BE_x\tau_D\asymp \delta(x),\quad x\in D.
\end{equation}
By the above (upper) estimate, $\langle m,R^D|\mu|\rangle\le c\langle|\mu|,\delta\rangle$,
so $\mu\in\RRr(D)$ if $\|\mu\|_{\delta,TV}<\infty$. Similarly, $f(\cdot,y)\in \RRr(D)$ if $f(\cdot,y)\in L^1_\delta(D;m)$.
In particular $f(\cdot,y)\in qL^1(D;m)$ in that case (see the remark following (\ref{eq6.6})).
Also note that if $g\in L^1(\partial D;\sigma)$,
then $P_D|g|(x)<\infty$ for every $x\in D$. Therefore, if (H1) is satisfied and
\begin{equation}
\label{eq7.10}
f(\cdot,y)\in L^1_\delta(D;m),\,\, y\in\BR,\qquad \mu\in\MMr_{\delta}(D),
\qquad g\in L^1(\partial D;\sigma),
\end{equation}
then
by Theorem \ref{th6.11} and Theorem \ref{th6.7} there exists a  unique solution $u$ of (\ref{eq1.1})
Let $\gamma=|f(\cdot,0)|\cdot m+|\mu|$. By  (\ref{eq7.6}) we have $\langle m,R^D\gamma\rangle =\langle\gamma,R^D1\rangle\le c\langle\gamma,\delta\rangle$ and $\langle m,R^D|f(\cdot,u)|\rangle\ge c^{-1}\langle|f(\cdot,u)|\cdot m,\delta\rangle$
for some $c\ge1$.
Moreover, by (\ref{eq7.5}), (\ref{eq7.6}) and Fubini's theorem,
\[
\|P_D(|g|)\|_{L^1(D;m)}\le c\int_{\partial D} \Big(\int_D|x-y|^{1-d}\,dx\Big)|g(y)|\,\sigma(dy).
\]
Therefore from (\ref{eq7.1}) it follows that $u\in L^1(D;m)$ and there is $C>0$ such that
\[
\|u\|_{L^1(D;m)}+\|f(\cdot,u)\|_{L^1_\delta(D;m)}\le C(\|f(\cdot,0)\|_{L^1_\delta(D;m)}+
\|\mu\|_{\MMr_{\delta}(D)}+\|g\|_{L^1(\partial D;\sigma)}).
\]
This means that if (H1) and (\ref{eq7.10}) is satisfied, then $u$ is a weak
solution in the sense defined in \cite[Section 1.2]{MV} and the estimate \cite[(2.1.8)]{MV}
holds true. For another proof of the existence
and uniqueness of weak solution of (\ref{eq1.1}) in case $L=\Delta$
and $\mu\in L^1(D;\delta\,dx)$
we refer the reader to \cite[Proposition 2.1.2]{MV}.
\end{example}

\begin{example}
\label{ex7.7}
Assume that $D$ is of class $C^{1,1}$, i.e. for every $y\in\partial D$ there is $r>0$
such that $B(x,r)\cap\partial D$ is the graph of a function whose first derivatives are
Lipschitz. We denote by $m$ or $dx$ the Lebesgue measure on $D$. Let  $L=-(-\Delta)^{\alpha/2}$
(see Example \ref{ex2.8}). In \cite[Proposition 4.9]{Ku} (see also \cite{CS}) it is proved that then
\begin{equation}
\label{eq7.7}
R^D1(x)=\BE_x\tau_D\asymp \delta^{\alpha/2}(x),\quad x\in D.
\end{equation}
In much the same way as in  Example \ref{ex7.6} from the upper estimate in (\ref{eq7.7})
we infer that if $f(\cdot,y)\in L^1_{\delta^{\alpha/2}}(D;m)$, $\|\mu\|_{\delta^{\alpha/2},TV}<\infty$,
then (H2) is satisfied and $\mu\in\RRr(D)$.
By \cite{Ku} (or \cite[Corollary 1.3]{CS}),
\begin{equation}
\label{eq7.8}
G_D(x,y)\asymp \min\Big\{\frac{1}{|x-y|^{-\alpha+d}}\,,
\frac{\delta^{\alpha/2}(x)\delta^{\alpha/2}(y)}{|x-y|^d}\Big\},\quad x,y\in D,
\end{equation}
and by \cite[Theorem 1.5]{CS},
\[
p_D(x,y)\asymp \frac{\delta^{\alpha/2}(x)}{\delta^{\alpha/2}(y)(1+\delta^{\alpha/2}(y))}\cdot
\frac{1}{|x-y|^d},\quad x\in D,\,y\in\BR^d\setminus\bar D,
\]
where $p_D$ is the density of the Poisson kernel, or, equivalently, $p_D$
is given by \eqref{eq.pkac} with $j(r)=c_{d,\alpha}r^{-\alpha-d}$
(see \cite[Theorem 1.4]{CS}).
Since $\delta(y)\le|x-y|$, it follows in particular that for some $c>0$,
\begin{equation}
\label{eq7.9}
p_D(x,y)\le c\delta^{\alpha/2}(x)\cdot\min\{\delta^{-\alpha/2}(y),\delta^{-\alpha-d}(y)\},
\quad x\in D,\,y\in\BR^d\setminus\bar D.
\end{equation}
We see that  $P_D|g|(x)<\infty$ for every $x\in D$ if $g\in L^1_{\delta^{-\alpha/2}}(D^c;m)$.
In particular (H3) is satisfied.  Thus, if (H1) is satisfied
and
\[
f(\cdot,y)\in L^1_{\delta^{\alpha/2}}(D;m),\,\,y\in\BR,\qquad \mu\in\MMr_{\delta^{\alpha/2}}(D),
\quad g\in L^1_{\delta^{-\alpha/2}}(D^c;m),
\]
then there exists a unique solution of (\ref{eq1.1}).
Furthermore, by (\ref{eq7.8}), for any nonnegative measure $\gamma$ on $D$ we have
\[
\int_D\Big(\int_DG_D(x,y)\,\gamma(dy)\Big)\,dx\le C\int_D\delta^{\alpha/2}(y)\,\gamma(dy),
\]
whereas by  (\ref{eq7.9}),
\[
\int_D\Big(\int_{D^c}|g(y)|p_D(x,y)\,dy\Big)dx
\le C\int_{D^c}|g(y)|\min\{\delta^{-\alpha/2}(y),\delta^{-\alpha-d}(y)\}\,dy.
\]
Observe also that by the lower estimate in (\ref{eq7.7}) there is $c>0$
such that  $\langle m,R^D|f(\cdot,u)|\ge c\langle |f(\cdot,u)|\cdot m,\delta^{\alpha/2}\rangle$.
Therefore, if
\begin{equation}
\label{eq7.12}
f(\cdot,y)\in L^1_{\delta^{\alpha/2}}(D;m),\,\,y\in\BR,\qquad \mu\in\MMr_{\delta^{\alpha/2}}(D)
\end{equation}
and
\begin{equation}
\label{eq7.13}
g\in L^1_{\delta^{-\alpha/2}\wedge \delta^{-\alpha-d}}(D^c;m)
\end{equation}
(cf. condition (11) in \cite{A}) then from (\ref{eq7.1}) it follows that $u\in L^1(D;dx)$ and
\begin{align*}
\|u\|_{L^1(D;m)}+\|f(\cdot,u)\|_{L^1_{\delta^{\alpha/2}}(D;m)}&\le C(\|f(\cdot,0)\|_{L^1_{\delta^{\alpha/2}}(D;m)}
+\|\mu\|_{\MMr_{\delta^{\alpha/2}}(D)} \\
&\quad+\|g\|_{L^1_{\delta^{-\alpha/2}\wedge \delta^{-\alpha-d}}(D^c;m)})
\end{align*}
for some $C>0$. This means that if (H1) and (\ref{eq7.12}), (\ref{eq7.13}) are satisfied, then $u$
is a weak solution, in the sense of \cite[Definition 1.3]{A}, of the problem
\[
-(-\Delta)^{-\alpha/2}=f(\cdot,u)+\mu\quad\mbox{in }D,
\quad u=g\quad\mbox{in }\BR^d\setminus\bar D,\quad Eu=0\quad\mbox{ on }\partial D,
\]
where $E$ is the boundary trace operator defined in \cite{A}.
\end{example}

\section{Weak and variational solutions}
\label{sec7b}

\begin{definition}[Variational solutions]
Let $g\in F$, $\mu\in F^*$. We say that $u$ is a variational  solution  of \eqref{eq1.1} if
$u\in F$, $u-g\in F(D)$, $f(\cdot,u)\in F^*$, and
for any $\eta\in F(D)$,
\[
\EE(u,\eta)=\int_D \eta f(\cdot,u)\,dm+\int_D \eta\,d\mu.
\]
\end{definition}

\begin{proposition}
\label{ex.prop.6.3}
Assume that $g\in F, \mu\in F^*$.
Let $u$ be a variational  solution  of \eqref{eq1.1}.
Then $u$ is a  solution of \eqref{eq1.1}.  Conversely, if $g\in F$, $\mu\in F^*$
and $u\in F$ is a solution
of \eqref{eq1.1} with  $f(\cdot,u)\in F^*$, then $u$ is a variational solution of \eqref{eq1.1}.
\end{proposition}
\begin{proof}
First note that for $g\in F$ we have
\begin{equation}
\label{eq5.gharm}
\EE(P_Dg,\eta)=-\EE(\Pi_Dg,\eta)+\EE(g,\eta)=0,\quad \eta\in F(D).
\end{equation}
Suppose  that $u$ is a  variational solution of \eqref{eq1.1}. Let $V_n:= D$, $n\ge 1$. With this  $D$-total family  all conditions
of Definition \ref{def6.1} are trivially satisfied (see Remark \ref{rem.equiv}).

Now, suppose that $u$ is a solution of \eqref{eq1.1} and $u,g\in F$, $f(\cdot,u),\mu\in F^*$.
By Theorem \ref{th6.7}, we may take $V_n=D$, $n \ge 1$. Consequently, $u$ is a variational solution of \eqref{eq1.1}.
\end{proof}

In case $L$ is purely nonlocal and the form associated with it admits the form
\[
\EE(u,v)=\int_{\mathbb R^d}\int_{\mathbb R^d}(u(x)-u(y))(v(x)-v(y))j(x,y)\,dx\,dy
\]
it is  natural, in the context of weak solutions, to use the notion of the space $V^D$
consisting of functions $u\in  L^2_{loc}(\mathbb R^d;m)$ such that
\[
\|u\|^2_{V^D}=2\int_{D\times \mathbb R^d}(u(x)-u(y))^2 j(x,y)\,dx\,dy<\infty.
\]
We let
\[
V^D(u,v)= 2\int_{D\times \mathbb R^d}(u(x)-u(y))(v(x)-v(y)) j(x,y)\,dx\,dy,\quad u,v\in V^D.
\]
This space  is widely used in the literature (see, e.g., \cite{DRV,DK,FKV,MSW}).
For simplicity of the presentation, we assume that $F\subset L^2(\mathbb R^d)$.
Clearly, $F(D)\subset V^D$, so $(V^D)^*\subset F^*(D)$.

We can now  introduce the notion of
weak solutions of \eqref{eq6.1} (see \cite{FKV}).

\begin{definition}[Weak solutions]
Let $\mu\in F^*(D)$ and $P_Dg\in V^D$. We say that $u\in V^D$ is a {\em weak solution}
to the problem
\begin{equation}
\label{eq6.1}
-Lu=\mu\quad\text{in }D,\qquad u=g\quad\text{on }\partial_{\chi}D,
\end{equation}
if $u=g$ q.e. on $\partial_{\chi}D$ and for any $\eta\in F(D)$
\begin{equation}
\label{eq1.1wle}
V^D(u,\eta)=\int_D\eta\,d\mu.
\end{equation}
\end{definition}

Our aim is  to apply the general existence result of Theorem \ref{th6.11}
to get, as an easy  corollary,   the existence result  for weak solutions of \eqref{eq1.1}.
Let us consider the following  conditions:
\begin{enumerate}
\item[(VD)] $F$ is dense in $V^D$.

\item[(VU)] If $u\in V^D$ and $V^D(u,u)=0$, then $u=0\,\,m$-a.e. in $\mathbb R^d$.
\end{enumerate}

\begin{theorem}
\label{th.srv}
Let $\mu\in F^*(D)$, $P_D|g|<\infty$ q.e. in $D$. Let $u$ be a solution of \eqref{eq1.1}
with $f=0$.
\begin{enumerate}
\item[\rm(i)] If  $P_Dg\in V^D$, then
 $u\in V^D$ and
\begin{equation}
\label{eq6.2wlq}
\EE(u-P_Dg,\eta)=\langle  \mu,\eta\rangle,\quad \eta\in F(D).
\end{equation}
Moreover,
\begin{equation}
\label{eq1.1wlein}
\|u\|_{V^D}\le \|P_Dg\|_{V^D}+\|\mu\|_{F^*(D)}.
\end{equation}
\item[\rm(ii)]  Assume \mbox{\rm(VD)}. If $P_Dg\in V^D$, then $u$ is a weak solution of \eqref{eq6.1}.
Furthermore, if $g\in V^D$, then $P_Dg\in V^D$ and
\begin{equation}
\label{eq1.1wlein12}
\|u\|_{V^D}\le \|g\|_{V^D}+\|\mu\|_{F^*(D)}.
\end{equation}
\end{enumerate}
\end{theorem}
\begin{proof}
(i) By Definition \ref{def6.1}, there exists a $D$-total family $(V_n)$
such that
\begin{equation}
\label{eq6.2wl}
\EE(\Pi_{V_n}(u),\eta)=\langle \mathbf1_{V_n}\cdot \mu,\eta\rangle
\end{equation}
for any $\eta\in F(V_n)$. Putting $\eta=\Pi_{V_n}(u)$ we get
\begin{equation}
\label{eq6.2wly}
\EE(\Pi_{V_n}(u),\Pi_{V_n}(u))\le \| \mu\|_{F^*(D)}.
\end{equation}
By Definition \ref{def6.1}(c), $\Pi_{V_n}(u)\to u-P_Dg$, which together with \eqref{eq6.2wly}
yields $u-P_Dg\in F(D)\subset  V^D$. Consequently, $u\in V^D$. By \eqref{eq6.2wly}, up to a subsequence,
$\Pi_{V_n}(u)\to u-P_Dg$ weakly in $F(D)$. Therefore letting $n\to\infty$ in \eqref{eq6.2wl} gives
\eqref{eq6.2wlq} for any $\eta\in \bigcup_{n\ge 1}F(V_n)$. By Proposition \ref{prop3.2}, \eqref{eq6.2wlq} holds for any $\eta\in F(D)$.
Therefore putting $\eta=u-P_Dg$ we obtain
\begin{equation}
\label{eq5.9}
\|u-P_Dg\|^2_{F(D)}=\langle  \mu,u-P_Dg\rangle\le\|u-P_Dg\|_{F(D)}\|\mu\|_{F^*(D)}.
\end{equation}
Observe now that by  Definition \ref{def6.1}(b), $\|u-P_Dg\|_{F(D)}=\|u-P_Dg\|_{V^D}$, which together with (\ref{eq5.9}) yields \eqref{eq1.1wlein}. This finishes the proof of (i). As for (ii), we observe that for $\eta\in F(D)$ and $g\in F$ we have
\[
V^D(P_Dg,\eta)=\EE(P_Dg,\eta)=0.
\]
Hence, by (VD), $V^D(P_Dg,\eta)=0$ for $g\in V^D$ and $\eta\in F(D)$.
This when combined with \eqref{eq6.2wlq} implies that $u$ is a weak solution of \eqref{eq6.1}.
Similarly, for $g\in F$ we have
\begin{align*}
&V^D(P_Dg,P_Dg)+\int_{D^c\times D^c}(P_Dg(x)-P_Dg(y))^2j(x,y)\,dx\,dy\\
&\qquad=\EE(P_Dg,P_Dg)
\le \EE(g,g)=V^D(g,g)+\int_{D^c\times D^c}(g(x)-g(y))^2j(x,y)\,dx\,dy.
\end{align*}
Since $P_Dg=g$ q.e. on $D^c$, we see that
\begin{align*}
V^D(P_Dg,P_Dg)\le V^D(g,g).
\end{align*}
From this and (VD) one easily concludes that $\|P_Dg\|_{V^D}\le\|g\|_{V^D}$ for $g\in V^D$,
which together with \eqref{eq1.1wlein} gives \eqref{eq1.1wlein12}.
\end{proof}

\begin{remark}
Observe that conditions (VU) and (VD) guarantee  uniqueness of weak solutions of \eqref{eq6.1}.
Therefore, under (VD) and (VU) we may use Theorem \ref{th.srv} and Theorem \ref{prop6.6}
to obtain a comparison result for weak solutions of \eqref{eq1.1}, i.e. for functions $u\in V^D$
such that $f(\cdot,u)\in F^*(D)$ and \eqref{eq1.1wle} holds with $\mu$ replaced by $f(\cdot,u)\cdot m+\mu$.
\end{remark}

\begin{remark}
Suppose that $g\in\mathcal B(D^c)$. In view of Theorem \ref{th.srv}, under (VD)
there exists a weak solution of \eqref{eq6.1} provided that $P_Dg\in V^D$.
One of the conditions guaranteeing that $P_Dg\in V^D$ is $g\in V^D$.
In general, this is not a necessary condition.
In \cite{BGPR}, under additional conditions on $j$, the authors give a necessary and sufficient
condition for  $g$ to satisfy $P_Dg\in V^D$. This  condition is expressed in terms
of appropriate Douglas integrals of $g$ over $D^c\times D^c$.
\end{remark}

\section{Very weak solutions}
\label{sec8}

We fix a bounded open set $D\subset E$, $\mu\in\MMr_{0,b}(D)$,
and $g\in\mathcal B(E)$ such that $P^D|g|\in L^1(E;m)$.
Throughout this section,
we assume that $R^D1$ is bounded. We set
\[
\mathfrak D_{[b]}(L)=\{\eta\in \mathfrak D(L)\cap \mathscr B_b(E): L\eta\in L^\infty(E;m)\}.
\]

In the present section, we focus on very weak solutions to the  problem (\ref{eq6.1}).

\begin{definition}
Let $\mathcal C$ be a  subset of $\mathfrak D_{[b]}(L)\cap F(D)$.
We say  that $u\in L^1(E;m)$ is a $\mathcal C$-{\em very weak solution}
of \eqref{eq6.1} if
\begin{equation}
\label{eq8.2cc}
-\int_E u\,L\eta\,dm =\int_D\eta\,d\mu,\quad  \eta\in \mathcal C,\qquad u=g\quad\text{q.e. on }\partial_{\chi}D.
\end{equation}
\end{definition}

\begin{lemma}
\label{lm8.1}
Suppose that $f\in \mathcal B_b(E)$ and $\eta\in \mathfrak D_{[b]}(L)\cap F(D)$.
Then
\begin{equation}
\label{eq8.3}
\int_E P_D(f)\,L\eta\,dm =0.
\end{equation}
\end{lemma}
\begin{proof}
Let
\[
\mathcal H=\{f\in\mathcal B_b(E): \eqref{eq8.3} \,\,\text{holds for any}\,\,\eta\in  \mathfrak D_{[b]}(L)\cap F(D) \}.
\]
Observe that $\mathcal H$ is a linear space and for any $(f_n)\subset \mathcal H$
such that $0\le f_n\le f_{n+1},\, n\ge 1$ a.e. we have $f:=\limsup_{n\to\infty}f_n\in \mathcal H$
provided $f$ is bounded. Next, observe that $F\cap \mathcal B_b(E)\subset \mathcal H$. Indeed, for any $f\in F$
and $\eta\in \mathfrak D(L)\cap F(D)$,
\[
\int_E P_D(f)\,(-L\eta)\,dm = \EE(P_D(f),\eta)=0.
\]
This implies in particular that $1\in \mathcal H$ ($1=\lim_{n\to \infty}e_{V_n}$ for any $E$-total family $(V_n)$
consisting of relatively compact open sets).
By \cite[Corollary 1.5.1]{FOT},  $F\cap C_b(E)$ is a multiplicative space, i.e.
$\eta_1\eta_2\in F\cap C_b(E)$ for all $\eta_1,\eta_2\in F\cap C_b(E)$.
Using the fact that $F\cap C_0(E)$
is dense in  $C_0(E)$ and the  monotone class theorem (see \cite[Theorem I.8]{P}) we get the desired result.
\end{proof}

\begin{theorem}
If $u$ is a solution of \eqref{eq1.1} with $f\equiv 0$ and  $u\in L^1(E;m)$, then it is a $\mathcal C$-very weak solution of \eqref{eq6.1} with $\mathcal C=\mathfrak D_{[b]}(L)\cap F(D)$.
\end{theorem}
\begin{proof}
Let $u$ be a solution of \eqref{eq6.1} and  $u\in L^1(E;m)$.
By Theorem \ref{th6.7}, $u=P_Dg+R^D\mu$ q.e.
By the assumptions we made, $P_Dg\in L^1(E;m)$.
Let $\eta\in \mathcal C$. Then
\[
\int_Du(-L\eta)\,dm=\int_DP_Dg(-L\eta)\,dm+\int_DR^D\mu(-L\eta)\,dm.
\]
The second term on the right-hand side of the above equality equals $\int_D\eta\,d\mu$.
Indeed,
\[
\int_DR^D\mu(-L\eta)\,dm=\int_DR^D(-L\eta)\,d\mu=\int_D \eta\,d\mu.
\]
Observe that $u=P_Dg=g$ q.e. in $D^c$. Hence
\[
\int_Eu(-L\eta)\,dm=\int_EP_Dg(-L\eta)\,dm+\int_D\eta\,d\mu.
\]
By Lemma \ref{lm8.1},
\[
\int_Eu(-L\eta)\,dm=\int_D\eta\,d\mu,
\]
which gives the result.
\end{proof}


\appendix

\section{Properties of orthogonal projections}

For  $u:E\to\bar\BR$ and $A\subset E$ we set
\begin{equation}
\label{eq.a1}
{\qesssup}_A u=\inf\{M\in[0,\infty]:u\le M\mbox{ q.e. in }A\}.
\end{equation}
Before proceeding, we note that for any $u\in F$ and constant $a>0$ the functions
$(u-a)^+$ and $-u^{-}$ are normal contractions of $u$, so $(u-a)^+, -u^{-}\in F$ by \cite[Corollary 1.5.1]{FOT}.
Moreover, for $u\in F$
\begin{equation}
\label{eq3.1}
\EE((u-a)^+,(u-a)^+)\le \EE(u,(u-a)^+), \quad \EE(-u^{-},-u^{-})\le\EE(u,-u^{-}).
\end{equation}
For $u\in\mathfrak D(\EE)$ the first inequality above follows from \cite[Theorem I.4.4]{MR} and the fact that $(u-a)^+=u-u\wedge a$.
The second inequality follows from the fact that $\EE(u^+,u^-)\le0$ (see \cite[p. 33]{MR}).
The case $u\in F$ follows by approximation.

\begin{proposition}
\label{prop3.2}
Let $g\in F$, and  $V\in\OO_q$. Then
\begin{enumerate}
\item[{\rm (i)}] $\qesssup_E|h_V(g)|\le \qesssup_{V^c}|g|$.

\item[{\rm (ii)}] $h_V(g)\ge 0$ q.e. provided $g\ge 0$ q.e.

\item[{\rm(iii)}] $\overline{\bigcup_{n\ge 1} F(V_n)}^F=F(V)$ for any sequence $(V_n)\subset\OO_q$ such that
    $V_n\uparrow V$ q.e.

\item[{\rm (iv)}] $\Pi_{V_n}(g)\rightarrow \Pi_V(g)$ in $F(V)$ for any $(V_n)$ as above.
\end{enumerate}
\end{proposition}
\begin{proof}
Set $a= \qesssup_{V^c}|g|$. Then $|g|\le a$ q.e. on $V^c$,  so by (\ref{eq2.14}),
$|h_V(g)|\le a$ q.e. on $V^c$. Since $(h_V(g)-a)^+\in F$, we see that
$(h_V(g)-a)^+\in F(V)$. By this, (\ref{eq3.1}) and the fact that $h_V(g)\in F(V)^{\bot}$ we get
\[
\EE((h_V(g)-a)^+,(h_V(g)-a)^+)\le \EE(h_V(g),(h_V(g)-a)^+)=0.
\]
This
implies that $h_V(g)\le a$ $m$-a.e. on $E$ and hence q.e. since $h_V(g)$ is quasi continuous.
In the same manner we can see that $h_V(g)\ge- a$ q.e., which proves (i).
If $g\ge0$ q.e., then $h_V(g)\ge0$ q.e. on $V^c$ by (\ref{eq2.14}). We also know that $-(h_V(g))^-\in F$.
Hence $-(h_V(g))^-\in F(V)$. On the other hand, $h_V(g)\in F(V)^{\bot}$. Therefore using (\ref{eq3.1}) we get
\[
\EE((h_V(g))^-,(h_V(g))^-)\le \EE(-(h_V(g))^-,-(h_V(g))^-)\le \EE(h_V(g),-(h_V(g))^-)=0,
\]
which together with quasi continuity of $h_V(g)$ implies (ii).
Let $f\in L^2(E;m)$. Observe that $\sup_{n\ge 1}\|R^{V_n}_\alpha f\|_F<\infty$, so
up to subsequence, $(R^{V_n}_\alpha f)$ converges weakly in $F$.
By \cite[Theorem 4.1]{Simon}, $\EE^{V_n}\to \EE^V$
in the strong resolvent sense, i.e. $R^{V_n}_\alpha f\to R^V_\alpha f$ in $L^2(E;m)$
for any $f\in L^2(E;m)$, which when combined with the weak convergence
of $(R^{V_n}_\alpha f)$  implies (iii) since $R^V_\alpha(L^2(E;m))$ is dense in $F(V)$
(see \cite[Lemma 1.3.3, Theorem 1.5.2(iii), Theorem 1.5.3($\gamma$)]{FOT}).
By the well known property  of the orthogonal projection,
\[
\EE(\Pi_{V_n}(g),\Pi_{V_n}(g))\le \EE(g,g),\quad  n\ge 1.
\]
Hence,  up to a subsequence, $\Pi_{V_n}(g)\to w$ weakly in $F$ for some $w\in \overline{\bigcup_{n\ge 1} F(V_n)}=F(V)$.
Consequently, for any $\eta\in  F(V_k)$,
\[
0=\EE(g- \Pi_{V_{n\vee k}}(g),\eta)\to \EE(g- w,\eta)\quad\mbox{as } n\to \infty.
\]
From this and (iii), $\EE(g- w,\eta)=0, \eta \in F(V)$, which  implies that $w=\Pi_V(g)$.
Consequently,  $\Pi_{V_n}(g)\rightarrow\Pi_V(g)$ weakly in $F$.
By this and the fact that $\Pi_{V_n}$, $\Pi_V$ are orthogonal projections we also have
$\EE(\Pi_{V_n}(g),\Pi_{V_n}(g))=\EE(\Pi_{V_n}(g),g)\rightarrow\EE(\Pi_V(g),g)
=\EE(\Pi_V(g),\Pi_V(g))$.
Using this and  the weak convergence of $\Pi_{V_n}(g)$ to $\Pi_V(g)$ again we get (iv).
\end{proof}

\section{Integrability properties of auxillary pocesses}
\label{app}

In what follows $\tau$ denotes a stopping time.

\begin{lemma}
\label{lem4.2}
\begin{enumerate}[\rm(i)]

\item  Suppose that $u\in\BB(E)$ and there exists a nonnegative measure $\nu\in\RRr(D)$
such that $|u|\le R^D\nu$ q.e. Then the family $\{u(X_\tau): \tau\le \tau_D\}$
is uniformly integrable  under the measure $\BP_x$ for q.e. $x\in D$.
\item If $g\in\mathcal B^n(E)$, $V\in \mathcal O_q$ and $P_V(|g|)<\infty$ q.e. in $V$,
then $\{P_V(g)(X_{\tau}):\tau\le\tau_V\}$ is uniformly integrable  under the measure $\BP_x$ for q.e. $x\in D$.
\end{enumerate}
\end{lemma}
\begin{proof}
(i) By \eqref{eq2.12}, $R^D\nu(x)=\mathbb E_xA^\nu_{\tau_{D}}$ q.e. in $D$.
By the strong Markov property and additivity of $A^\nu$, for q.e. $x\in D$ and for any stopping time
$\tau\le\tau_D$,
\[
|u(X_\tau)|\le \mathbb E_x(A^\nu_{\tau_D}-A^\nu_\tau|\FF_\tau)\le  \mathbb E_x(A^\nu_{\tau_D}|\FF_\tau)\quad \BP_x\text{-a.s.}
\]
This readily yields (i).
\\
(ii) By the strong Markov property, for q.e. $x\in D$ and for any stopping time
$\tau\le\tau_V$,
\[
P_V(g)(X_{\tau})=\mathbb E_x(g(X_{\tau_V})|\FF_{\tau})\quad \BP_x\text{-a.s.}
\]
From this one easily deduces the result.
\end{proof}

\begin{lemma}
\label{lm.itoinc}
Let $x\in E$,   $\tau$ be a stopping time and $A$ be a continuous $\mathbb F$-adapted increasing process such that  $A_0=0$ under the measure $\mathbb P_x$.
\begin{enumerate}
\item[\rm(i)] If $\mathbb E_xA^2_\tau<\infty$, then
\[
\mathbb E_x(A_{\tau})^2=2\mathbb E_x\int_0^{\tau}\mathbb E_x(A_{\tau}-A_t|\FF_t)\,dA_t.
\]
\item[\rm(ii)] If $\mathbb E_xA_\tau<\infty$ and there exists $c>0$ such that $\mathbb E_x(A_\tau-A_t|\FF_t)\le c$,
$\mathbb P_x$-a.s. for any $t\ge 0$, then $\mathbb E_xA^2_\tau<\infty$.
\end{enumerate}
\end{lemma}
\begin{proof}
(i) Let $N^x$ denote a c\`adl\`ag version of the
martingale ${\mathbb{E}}_x (A_{\tau}|\FF_t)$, $t\ge0$.
By the Doob $L^2$-inequality, $\mathbb E_x\sup_{t\le\tau}|N^x_t|^2\le 4\mathbb E_xA^2_\tau<\infty$.
Consequently, by the Burkholder--Davis--Gundy inequality, $M^x_t=\int_0^{t\wedge\tau}A_s\,dN^x_s$ is a martingale.
Integrating by parts we obtain
\begin{align*}
2{\mathbb{E}}_x\int^{\tau}_0{\mathbb{E}}_x(A_{\tau}
-A_t|{\mathcal{F}}_t)\,dA_t
&=2{\mathbb{E}}_x\int^{\tau}_0N^x_t\,dA_t
-2{\mathbb{E}}_x\int^{\tau}_0A_t\,dA_t\\
&=2{\mathbb{E}}_x\Big(N^x_{\tau}A_{\tau}-\int^{\tau}_0A_t\,dN^x_t
-\frac12(A_{\tau})^2\Big)={\mathbb{E}}_x(A_{\tau})^2.
\end{align*}
This completes the proof of (i). For (ii), let $\tau_n=\inf\{t\ge 0: A_t\ge n\}$.
By (i) and the assumptions,
\[
\mathbb E_x(A_{\tau_n\wedge \tau})^2=2\mathbb E_x\int_0^{\tau_n\wedge\tau}\mathbb E_x(A_{\tau_n\wedge\tau}-A_t|\FF_t)\,dA_t\le 2c \mathbb E_xA_{\tau}.
\]
Applying Fatou's lemma gives (ii).
\end{proof}

\subsection*{Acknowledgements}
The first author was supported  by Polish National Science Centre
under grant no. 2017/25/B/ST1/00878.



\end{document}